\documentclass[a4paper]{amsart}

\usepackage{amsfonts}
\usepackage{amssymb}
\usepackage{amsmath}
\usepackage{hyperref}
\usepackage{mathrsfs}
\usepackage{centernot}
\usepackage{mathdots}
\usepackage{stmaryrd}
\usepackage[all]{xy}

\newtheorem{thm}{Theorem}[section]

\newtheorem{prp}[thm]{Proposition}
\newtheorem{cor}[thm]{Corollary}
\newtheorem{dfn}[thm]{Definition}

\newtheoremstyle{roman} 
    {8.0pt plus 2.0pt minus 4.0pt}                    
    {8.0pt plus 2.0pt minus 4.0pt}                    
    {\normalfont}                
    {}                           
    {\bfseries}                  
    {.}                          
    {5pt plus 1pt minus 1pt}     
    {}  

\theoremstyle{roman}

\newtheorem{example}[thm]{Example}
\newtheorem{remark}[thm]{Remark}

\theoremstyle{plain}

\newcommand{\rem}[1]{}

\newcommand{\C}{\mathbb{C}}

\newcommand{\N}{\mathbb{N}}
\newcommand{\Q}{\mathbb{Q}}
\newcommand{\R}{\mathbb{R}}
\newcommand{\Z}{\mathbb{Z}}

\newcommand{\frakCapital}{
\newcommand{\frakA}{{\mathfrak{A}}}
\newcommand{\frakB}{{\mathfrak{B}}}
\newcommand{\frakC}{{\mathfrak{C}}}
\newcommand{\frakD}{{\mathfrak{D}}}
\newcommand{\frakE}{{\mathfrak{E}}}
\newcommand{\frakF}{{\mathfrak{F}}}
\newcommand{\frakG}{{\mathfrak{G}}}
\newcommand{\frakH}{{\mathfrak{H}}}
\newcommand{\frakI}{{\mathfrak{I}}}
\newcommand{\frakJ}{{\mathfrak{J}}}
\newcommand{\frakK}{{\mathfrak{K}}}
\newcommand{\frakL}{{\mathfrak{L}}}
\newcommand{\frakM}{{\mathfrak{M}}}
\newcommand{\frakN}{{\mathfrak{N}}}
\newcommand{\frakO}{{\mathfrak{O}}}
\newcommand{\frakP}{{\mathfrak{P}}}
\newcommand{\frakQ}{{\mathfrak{Q}}}
\newcommand{\frakR}{{\mathfrak{R}}}
\newcommand{\frakS}{{\mathfrak{S}}}
\newcommand{\frakT}{{\mathfrak{T}}}
\newcommand{\frakU}{{\mathfrak{U}}}
\newcommand{\frakV}{{\mathfrak{V}}}
\newcommand{\frakW}{{\mathfrak{W}}}
\newcommand{\frakX}{{\mathfrak{X}}}
\newcommand{\frakY}{{\mathfrak{Y}}}
\newcommand{\frakZ}{{\mathfrak{Z}}}
}

\newcommand{\calCapital}{
\newcommand{\calA}{{\mathcal{A}}}
\newcommand{\calB}{{\mathcal{B}}}
\newcommand{\calC}{{\mathcal{C}}}
\newcommand{\calD}{{\mathcal{D}}}
\newcommand{\calE}{{\mathcal{E}}}
\newcommand{\calF}{{\mathcal{F}}}
\newcommand{\calG}{{\mathcal{G}}}
\newcommand{\calH}{{\mathcal{H}}}
\newcommand{\calI}{{\mathcal{I}}}
\newcommand{\calJ}{{\mathcal{J}}}
\newcommand{\calK}{{\mathcal{K}}}
\newcommand{\calL}{{\mathcal{L}}}
\newcommand{\calM}{{\mathcal{M}}}
\newcommand{\calN}{{\mathcal{N}}}
\newcommand{\calO}{{\mathcal{O}}}
\newcommand{\calP}{{\mathcal{P}}}
\newcommand{\calQ}{{\mathcal{Q}}}
\newcommand{\calR}{{\mathcal{R}}}
\newcommand{\calS}{{\mathcal{S}}}
\newcommand{\calT}{{\mathcal{T}}}
\newcommand{\calU}{{\mathcal{U}}}
\newcommand{\calV}{{\mathcal{V}}}
\newcommand{\calW}{{\mathcal{W}}}
\newcommand{\calX}{{\mathcal{X}}}
\newcommand{\calY}{{\mathcal{Y}}}
\newcommand{\calZ}{{\mathcal{Z}}}
}

\newcommand{\bbCapital}{
\newcommand{\bbA}{{\mathbb{A}}}
\newcommand{\bbB}{{\mathbb{B}}}
\newcommand{\bbC}{{\mathbb{C}}}
\newcommand{\bbD}{{\mathbb{D}}}
\newcommand{\bbE}{{\mathbb{E}}}
\newcommand{\bbF}{{\mathbb{F}}}
\newcommand{\bbG}{{\mathbb{G}}}
\newcommand{\bbH}{{\mathbb{H}}}
\newcommand{\bbI}{{\mathbb{I}}}
\newcommand{\bbJ}{{\mathbb{J}}}
\newcommand{\bbK}{{\mathbb{K}}}
\newcommand{\bbL}{{\mathbb{L}}}
\newcommand{\bbM}{{\mathbb{M}}}
\newcommand{\bbN}{{\mathbb{N}}}
\newcommand{\bbO}{{\mathbb{O}}}
\newcommand{\bbP}{{\mathbb{P}}}
\newcommand{\bbQ}{{\mathbb{Q}}}
\newcommand{\bbR}{{\mathbb{R}}}
\newcommand{\bbS}{{\mathbb{S}}}
\newcommand{\bbT}{{\mathbb{T}}}
\newcommand{\bbU}{{\mathbb{U}}}
\newcommand{\bbV}{{\mathbb{V}}}
\newcommand{\bbW}{{\mathbb{W}}}
\newcommand{\bbX}{{\mathbb{X}}}
\newcommand{\bbY}{{\mathbb{Y}}}
\newcommand{\bbZ}{{\mathbb{Z}}}
}

\newcommand{\bfCapital}{
\newcommand{\bfA}{{\mathbf{A}}}
\newcommand{\bfB}{{\mathbf{B}}}
\newcommand{\bfC}{{\mathbf{C}}}
\newcommand{\bfD}{{\mathbf{D}}}
\newcommand{\bfE}{{\mathbf{E}}}
\newcommand{\bfF}{{\mathbf{F}}}
\newcommand{\bfG}{{\mathbf{G}}}
\newcommand{\bfH}{{\mathbf{H}}}
\newcommand{\bfI}{{\mathbf{I}}}
\newcommand{\bfJ}{{\mathbf{J}}}
\newcommand{\bfK}{{\mathbf{K}}}
\newcommand{\bfL}{{\mathbf{L}}}
\newcommand{\bfM}{{\mathbf{M}}}
\newcommand{\bfN}{{\mathbf{N}}}
\newcommand{\bfO}{{\mathbf{O}}}
\newcommand{\bfP}{{\mathbf{P}}}
\newcommand{\bfQ}{{\mathbf{Q}}}
\newcommand{\bfR}{{\mathbf{R}}}
\newcommand{\bfS}{{\mathbf{S}}}
\newcommand{\bfT}{{\mathbf{T}}}
\newcommand{\bfU}{{\mathbf{U}}}
\newcommand{\bfV}{{\mathbf{V}}}
\newcommand{\bfW}{{\mathbf{W}}}
\newcommand{\bfX}{{\mathbf{X}}}
\newcommand{\bfY}{{\mathbf{Y}}}
\newcommand{\bfZ}{{\mathbf{Z}}}
}

\newcommand{\catCapital}{
\newcommand{\catA}{{\mathscr{A}}}
\newcommand{\catB}{{\mathscr{B}}}
\newcommand{\catC}{{\mathscr{C}}}
\newcommand{\catD}{{\mathscr{D}}}
\newcommand{\catE}{{\mathscr{E}}}
\newcommand{\catF}{{\mathscr{F}}}
\newcommand{\catG}{{\mathscr{G}}}
\newcommand{\catH}{{\mathscr{H}}}
\newcommand{\catI}{{\mathscr{I}}}
\newcommand{\catJ}{{\mathscr{J}}}
\newcommand{\catK}{{\mathscr{K}}}
\newcommand{\catL}{{\mathscr{L}}}
\newcommand{\catM}{{\mathscr{M}}}
\newcommand{\catN}{{\mathscr{N}}}
\newcommand{\catO}{{\mathscr{O}}}
\newcommand{\catP}{{\mathscr{P}}}
\newcommand{\catQ}{{\mathscr{Q}}}
\newcommand{\catR}{{\mathscr{R}}}
\newcommand{\catS}{{\mathscr{S}}}
\newcommand{\catT}{{\mathscr{T}}}
\newcommand{\catU}{{\mathscr{U}}}
\newcommand{\catV}{{\mathscr{V}}}
\newcommand{\catW}{{\mathscr{W}}}
\newcommand{\catX}{{\mathscr{X}}}
\newcommand{\catY}{{\mathscr{Y}}}
\newcommand{\catZ}{{\mathscr{Z}}}
}

\newcommand{\what}[1]{\widehat{#1}}

\newcommand{\veps}{\varepsilon}
\newcommand{\vphi}{\varphi}

\newcommand{\idealof}{\unlhd} 

\newcommand{\onto}{\twoheadrightarrow}

\newcommand{\suchthat}{\,:\,}
\newcommand{\where}{\,|\,}

\newcommand{\quo}[1]{\overline{#1}}

 \newcommand{\DDotsArr}[2]{
\begin{array}{ccc} {#1} &  &  \\
 & \ddots &  \\
 &  & {#2} \end{array}}

\newcommand{\Trings}[1]{\left< #1 \right>}

 \DeclareMathOperator{\Aut}{Aut}
\DeclareMathOperator{\Char}{char} %
\DeclareMathOperator{\End}{End} %
\DeclareMathOperator{\Hom}{Hom} %
\DeclareMathOperator{\id}{id} %
\DeclareMathOperator{\im}{im} %
\DeclareMathOperator{\Nrd}{Nrd} %
\newcommand{\op}{\mathrm{op}} %
\DeclareMathOperator{\Span}{span} %
\DeclareMathOperator{\Spec}{Spec} %
\DeclareMathOperator{\Stab}{Stab} %
\newcommand{\nGL}[2]{\mathrm{GL}_{#2}({#1})}
\newcommand{\nPGL}[2]{\mathrm{PGL}_{#2}({#1})}

\newcommand{\nMat}[2]{\mathrm{M}_{#2}(#1)}

\newcommand{\uGL}{{\mathbf{GL}}}
\newcommand{\uPGL}{{\mathbf{PGL}}}
\newcommand{\uSL}{{\mathbf{SL}}}

\newcommand{\units}[1]{{#1^\times}}

\usepackage[foot]{amsaddr}

\usepackage{xcolor}

\usepackage{dsfont} 



\calCapital %
\frakCapital %
\bbCapital %
\catCapital %
\bfCapital %

\newcommand{\Rep}[2][]{\mathrm{Rep}^{\mathrm{#1}}(#2)}
\newcommand{\Irr}[2][]{\mathrm{Irr}^{\mathrm{#1}}(#2)}

\DeclareMathOperator{\Ball}{B}
\DeclareMathOperator{\dist}{d}

\newcommand{\leftmod}{ {\setminus} }

\newcommand{\vrt}[1]{{#1}_{\mathrm{vrt}}}

\newcommand{\ori}[1]{#1_{\mathrm{ori}}}

\newcommand{\Alg}[3][]{A^{#1}_{#3}({#2})}
\newcommand{\tAlg}[3][]{\what{A}^{#1}_{#3}({#2})}

\newcommand{\cconj}[1]{\quo{#1}} 

\newcommand{\catSimp}{\mathtt{Cov}}

\newcommand{\catPHil}{\mathtt{pHil}}
\newcommand{\catSet}{\mathtt{Set}}

\newcommand{\udual}[1]{{\what{#1}}}

\newcommand{\norm}[1]{\|{#1}\|}

\newcommand{\abs}[1]{|{#1}|}
\newcommand{\Abs}[1]{\left|{#1}\right|}

\newcommand{\wc}{\prec}
\newcommand{\sphere}[2][1]{\mathbb{S}^{#1}({#2})}
\newcommand{\fsubseteq}{\subseteq_f}

\newcommand{\LL}[2][]{\mathrm{L}^{2}_{#1}({#2})}

\newcommand{\llf}{\tilde{\ell}^2}

\newcommand{\e}{\mathrm{e}}

\newcommand{\sfx}{{\mathsf{x}}}

\newcommand{\change}[2]{#2}

\newcommand{\bGL}{\uGL}
\newcommand{\bPGL}{\uPGL}

\numberwithin{equation}{section}

\title[The Ramanujan Property for Simplicial Complexes]{Highlights from ``The Ramanujan Property for Simplicial Complexes''}

\author{Uriya A.\ First$^*$}
\date{\today}
\address{$^*$Department of Mathematics, University of British Columbia}
\email{uriya.first@gmail.com}
\thanks{This research was supported by an ERC grant \#226135,
the Lady Davis Fellowship Trust, and the UBC Mathematics Department.}


\begin{document}

\maketitle



\begin{abstract}
    This paper brings the main definitions and results from ``The Ramanujan Property for Simplicial Complexes''.
    No proofs are given.
    
    Given a simplicial complex $\calX$ and a group $G$ acting on $\calX$, we define
    Ramanujan quotients $\calX$. For $G$ and $\calX$ suitably
    chosen this recovers \emph{Ramanujan $k$-regular graphs}
    and \emph{Ramanujan complexes} in the sense of Lubotzky, Samuels and Vishne. 
    Deep results in automorphic representations are used to give new examples of Ramanujan quotients when $\calX$
    is the affine building of an inner form of $\uGL_n$ over a local field of positive characteristic.
\end{abstract}




\bigskip

	This is an overview of the paper
	``The Ramanujan Property for Simplicial Complexes'' \cite{Fi15A}.
	We bring all main results and definitions (simplified at times), but give no proofs.
	Rather, precise references to the relevant
	places in \cite{Fi15A} are provided. 
	

\section*{Introduction}

	Let $X$ be a connected $k$-regular graph.
	Denote by $\lambda(X)$  the maximal absolute value of an eigenvalue of the adjacency matrix of $X$, excluding
    $k$ and $-k$. 
    Graphs for which $\lambda(X)$ is bounded away from $k$ are called \emph{expander graphs}. They enjoy many
    good combinatorial properties; see \cite{Lub12A} for a survey.
    
    The graph $X$ is called \emph{Ramanujan} if $\lambda(X)\leq 2\sqrt{k-1}$.
    This definition is motivated by the Alon--Boppana Theorem \cite{Nilli91}, stating that for any $\veps>0$, only finitely
    many non-isomorphic $k$-regular graphs satisfy the tighter bound $\lambda(X)\leq 2\sqrt{k-1}-\veps$.
    Furthermore,  the interval $[-2\sqrt{k-1},2\sqrt{k-1}]$ is the spectrum of the  $k$-regular
    tree (\cite[p.~252, Apx.~3]{Sunada88}, \cite[Th.~3]{Kest59}), which is the universal cover of any $k$-regular graph.
	Ramanujan graphs can therefore be thought of as having the smallest possible spectrum one can expect
    of an infinite family of graphs, or as  finite approximations of  the infinite $k$-regular tree.
    
    The spectral properties of
    Ramanujan graphs manifest in some desired combinatorial properties.
    For example, Ramanujan graphs are  supreme expanders, and  have a large chromatic number if not bipartite.
    Some known
    constructions
    have large girth as well.
    Constructing infinite families of non-isomorphic $k$-regular Ramanuajan graphs  is considered difficult.
    The first such families were introduced by Lubotzky, Phillips and Sarnak \cite{LubPhiSar88} and
    independently  by
    Margulis \cite{Marg88},
    assuming $k-1$ is prime. Morgenstern \cite{Morg94} has extended this to the case $k-1$ is a prime power.
    These works rely on deep results of Delinge \cite{Deligne74} and Drinfeld \cite{Drinfel88} concerning
    the Ramanujan--Petersson conjecture for $\bGL_2$.
    The existence of infinitely many $k$-regular bipartite Ramanujan graphs for arbitrary $k$ was later shown  by Marcus, Spielman and
    Srivastava \cite{MarSpiSri14} using different methods; the non-bipartite case remains open.

\medskip

	A high-dimensional generalization of Ramanujan graphs, called \emph{Ramanujan complexes},
    was suggested
    by Cartwright, Sol\'{e} and \.{Z}uk \cite{CarSolZuk03}, and later
    refined by Lubotzky, Samuels and Vishne \cite{LubSamVi05}
    (see \cite{JorLiv00}
    for another generalization of  Ramanujan graphs).
    These complexes are quotients of the affine Bruhat--Tits building of $\nPGL{F}{d}$,
    denoted  $\calB_d(F)$,
    where
    $F$ is a non-archimedean local field. The building $\calB_d(F)$ is a contractible
    simplicial complex of dimension $d-1$; its construction is recalled in \ref{subsec:building} below.
    The spectrum of a quotient
    of $\calB_d(F)$, i.e.\ a simplicial complex whose universal cover is $\calB_d(F)$,
    consists of the common spectrum of a certain
    family of $d-1$ linear operators associated with the quotient, called the \emph{Hecke operators}.
    According to Lubotzky, Samuels and Vishne  \cite{LubSamVi05}, a quotient of $\calB_d(F)$ is \emph{Ramanujan} if its spectrum,
    which is a subset of $\C^{d-1}$, is contained in the spectrum of the universal cover $\calB_d(F)$ together with a certain family
    of $d$ points in $\C^{d-1}$, called the \emph{trivial spectrum}.

    Li \cite[Thm.~4.1]{Li04} proved a theorem in the spirit of the Alon--Boppana Theorem for quotients of $\calB_d(F)$:
    If $\{X_n\}_{n\in\N}$ is a family of such quotients satisfying a mild assumption,
    then the closure of the union of the spectra of $\{X_n\}_{n\in\N}$ (in $\C^{d-1}$) contains
    the spectrum of the universal cover $\calB_d(F)$. Ramanujan complexes can therefore be thought of as having
    the smallest possible spectrum that can be expected of an
    infinite family, or  as  spectral approximations of the universal cover
    $\calB_d(F)$, similarly to Ramanujan graphs. When $d=2$, the complex $\calB_d(F)$ is a regular tree,
    and  Ramanujan complexes
    are just Ramanujan graphs in the previous sense.

    The
    existence of infinite families of Ramanujan complexes was shown by Lubotzky, Samuels and Vishne in \cite{LubSamVi05}
    (see also~\cite{LubSamVish05B}),
    using Lafforgue's proof of the Rama\-nu\-jan--Peters\-son conjecture for $\bGL_d$ in positive characteristic  \cite{Laff02}.\change{4}{\footnote{
    	The proof in \cite{LubSamVi05} assumed the global Jacquet--Langlands correspondence
    	for $\uGL_n$
    	in positive characteristic that was established later in \cite{BadulRoch14}.
    }
    Li \cite{Li04} has independently obtained very similar results using
    a special case of the conjecture established by Laumon, Rapoport and Stuhler \cite[Th.~14.12]{LaRaSt93}.\footnote{
    	The notion of Ramanujan complexes used in \cite{Li04} is slightly
    	weaker than the one used in \cite{LubSamVi05}, but the constructions of  \cite{Li04} are in fact Ramanujan
    	in the sense of \cite{LubSamVi05}.
	}}
    As in the case of graphs, Ramanujan complexes enjoy various good combinatorial properties:
    They have high chromatic number \cite[\S6]{EvraGolLub14}, good mixing properties \cite[\S4]{EvraGolLub14},
    they
    satisfy Gromov's \emph{geometrical expansion property} \cite{FoGrLaNaPa12} (see also \cite{Gro10},
    \cite{KauKazhLub14}), and the constructions of \cite{LubSamVi05} have high girth in addition \cite{LubMesh07}.

\medskip

	The Ramanujan property of quotients of $\calB_d(F)$ is measured with respect to the spectrum
    of the \emph{Hecke operators}. In a certain sense, to be made precise 
    below,
    these operators
    capture all spectral information in dimension $0$. Therefore, we regard the spectrum of Lubotzky, Samuels and Vishne
    as the \emph{$0$-dimensional spectrum}. However, one can associate other operators with a simplicial complex such
    that their spectrum
    affects combinatorial properties. For example, this is the case for the high-dimensional Laplacians; see
    for instance \cite{ParRosTes13}, \cite{Golubev13},
    \cite{GolPar14}, \cite{Par15}.
    Other examples are adjacency operators
    between  various types of facets.
    These operators are high-dimensional in nature and so
    their spectrum is a
    priori not determined by the spectrum 
    of the  Hecke operators.
    The work \cite{Fi15A}, which is summarized in this manuscript, treats  these and other high-dimensional operators, and
    constructs  examples of complexes which are Ramanujan relative to such operators.
    
\medskip

	In more detail, let $\calX$ be a simplicial complex and let $G$ be a group of automorphisms of $\calX$ satisfying certain mild assumptions.
    For example, one can take $\calX$ to be a  $k$-regular
    tree  $\calT_k$ and $G=\Aut(\calT_k)$, or $\calX=\calB_d(F)$ and $G=\nPGL{F}{d}$.
    Even more generally, $\calX$ can be an affine Bruhat-Tits building (see \cite{Buildings08AbramBrown}), and
    $G$ can be a group of automorphisms acting  on $\calX$ in a sufficiently transitive manner.
    We consider quotients of $\calX$ by  subgroups of $G$, called $G$-quotients for brevity,
    and associate several types of spectra with each of them.
    Among these spectra is the
    \emph{(non-oriented) $i$-dimensional spectrum}.
    When $\calX=\calT_k$
    and $G=\Aut(\calT_k)$, or when $\calX=\calB_d(F)$ and $G=\nPGL{F}{d}$, our $0$-dimensional spectrum coincides with
    the spectra of  quotients of regular graphs and quotients of $\calB_d(F)$  discussed earlier.
    
\medskip
    
    {\bf The main results} of \cite{Fi15A} are as follows:
    \begin{enumerate}
    	\item[(1)] 
    	If $\{X_n\}_{n\in\N}$ is a family of  $G$-quotients of $\calX$ satisfying a mild assumption, then
    	the closure of $\bigcup_{n\in\N}\Spec(X_n)$ contains $\Spec(\calX)$ (Theorem~\ref{TH:Alon-Boppana-I}).
    \end{enumerate}
    This generalizes Li's aforementioned theorem, and
    leads to a notion of \emph{Ramanujan $G$-quotients} of $\calX$:
    \begin{enumerate}
    	\item[($*$)] A $G$-quotient of $\calX$ is Ramanujan (relative to a particular type of spectrum)
    	if its spectrum is contained in the union of the spectrum of $\calX$ with the \emph{trivial spectrum}
    	(Definition~\ref{DF:Ramanujan-quotient}).
    \end{enumerate} 
    In analogy with Ramanujan graphs and
    Ramanujan complexes, Ramanujan $G$-quotients  have the smallest possible spectrum one can expect of an infinite family of $G$-quotients
    of $\calX$, or alternatively, they can be regarded as spectral approximations of the covering complex $\calX$.
    When $\calX$ is a $k$-regular tree (resp.\ $\calB_d(F)$), the quotients of $\calX$ which
    are \emph{Ramanujan in dimension $0$} are precisely the Ramanujan graphs (resp.\ Ramanujan complexes
    in the sense of \cite{LubSamVi05}).
    
    Next, we give a representation-theoretic criterion for being Ramanujan (Theorem~\ref{TH:Ramanujan-criterion}):
    Let  $\Gamma\leq G$ be a subgroup such that $\Gamma\leftmod\calX$
    is a finite simplicial complex and $\calX\to \Gamma\leftmod\calX$ is a cover map. 
    \begin{enumerate}
    	\item[(2)] Let $x_1,\dots,x_t$ be representatives for the orbits
    	of the action of $G$ on the $i$-dimensional cells in $\calX$, and write $K_n=\Stab_G(x_n)$.
    	Then $\calX$ is \emph{Ramanujan in dimension $i$}
    	if and only if every
    	irreducible unitary $G$-subrepresentation of $\LL{\Gamma\leftmod G}$ with $V^{K_1}+\dots+V^{K_t}\neq 0$
    	is 
    	\emph{tempered} or finite-dimensional.
    	\item[(2$'$)] $\Gamma\leftmod \calX$ is \emph{completely Ramanujan} if and only if every
    	irreducible unitary $G$-subrepresentation of $\LL{\Gamma\leftmod G}$ is 
    	\emph{tempered} or finite-dimensional.
	\end{enumerate} 
	The completely Ramanujan condition means that the complex is Ramanujan with respect to ``any type''
	of spectrum. In particular, it is {Ramanujan in all dimensions} and  Ramanujan relative to all
	high-dimensional Laplacians. However, it is also Ramanujan relative to other operators
	such as
	the adjacency operator of the graph obtained from $\Gamma\leftmod \calX$ by taking triples of vertices $(u_1,u_2,u_3)$
	with $\dist(u_1,u_2)=\dist(u_2,u_3)=1$ and saying that $(u_1,u_2,u_3)$ is adjacent to $(v_1,v_2,v_3)$
	if $u_2=v_1$ and $u_3=v_2$.\footnote{
		It should be pointed out that this does not mean that the regular graph  obtained in this manner
		is Ramanujan. Rather, the spectrum of its adjacency operator is contained in the union of
		the spectrum
		of the adjacency operator of the graph obtained from $\calX$ and the \emph{trivial spectrum}.
	}
	The result (2) can  be refined
	into a one-to-one correspondence between  the $i$-dimensional spectrum of $\Gamma\leftmod\calX$
    and a certain class of $G$-subrepresentations of $\LL{\Gamma\leftmod G}$ (Theorem~\ref{TH:spectrum-correspondence}).
	
    In case $\calX$ is the affine Bruhat--Tits building of a simple algebraic group $\bfG$ over the global field
    of $F$,
    $G=\bfG(F)$, and $\Gamma$ is an arithmetic cocompact lattice in $G$, our criterion can be restated
    in terms of automorphic representations of $\bfG$
    (Theorem~\ref{TH:automorphic-ramanujan}). The latter is used together with 
    deep results about automorphic representations (particularly \cite{Laff02} and \cite{BadulRoch14})
    to show:
    \begin{enumerate}
    	\item[(3)] Let $F$ be a non-archimedean local field with $\Char F>0$,
    	let $D$ be a central division algebra over $F$, let $d\geq 2$,
    	let $G=\nPGL{D}{d}:=\nGL{D}{d}/\units{F}$, and let $\calB_d(D)$
    	be the affine Bruhat--Tits building of $G$ (cf.\ \ref{subsec:building}). 
    	Then there exist 
    	infinite families of $G$-quotients of $\calB_d(D)$ which are
    	completely Ramanujan (Theorem~\ref{TH:ram-quo-exists}).
    \end{enumerate}
    Particular $G$-quotients of $\calB_d(D)$ which are completely Ramanujan
    are constructed in \ref{subsec:ramanujan-exist}.
    When $D=F$, our Ramanujan $G$-quotients are the Ramanujan
    complexes constructed by Lubotzky, Samuels and Vishne \cite{LubSamVi05}. Thus, the Ramanujan complexes
    of  \cite{LubSamVi05}, which are Ramanujan in dimension $0$ according to our setting,
    are in fact completely Ramanujan.
    When $d=2$, our construction gives rise to  Ramanujan graphs, which seem to be new when $D\neq F$.

    \medskip
    
    The paper is organized as follows:
    Section~\ref{sec:algebras} is a brief introduction to the spectral theory
    of $*$-algebras. It brings some auxiliary definitions and results used later
    in the text.
    Section~\ref{sec:simp-comp} recalls simplicial complexes and certain facts about $\ell$-groups acting on them.
    In Section~\ref{sec:spectrum}, we introduce our notion of  spectrum together
    with examples and supplementary results. 
    In Section~\ref{sec:opt-spec}, we present a generalization of Li's Theorem (reminiscent of the Alon--Boppana Theorem)
    and introduce the \emph{trivial spectrum}, which leads to the definition of \emph{Ramanujan quotients}.
    Section~\ref{sec:rep-thy} gives a 
    representation-theoretic criterion for a quotient of $\calX$ to be Ramanujan.
    Some consequence are  discussed.
    Finally, in Section~\ref{sec:ram-comp}, we recall the construction
    of the affine Bruhat--Tits building of $\nPGL{D}{d}$ and describe an infinite
    family of
    completely Ramanujan quotients of it, provided $\Char F>0$. 
    We also  explain how the problem is translated
    into a statement about automorphic representations.
    
\rem{
    \subsection*{Acknowledgements}

    We owe a debt of gratitude to Alex Lubotzky for presenting us with the theory of Ramanujan graphs
    and suggesting this research project. We are also in debt to Lior Silberman for many beneficial discussions.
    We further thank Anne-Marie Aubert, Alexandru Ioan Badulescu,
    David Kazhdan, Laurent Lafforgue and Dipendra Prasad for short-yet-crucial   correspondences,
    all related to automorphic representations. In addition, Amitay Kamber has given us several useful suggestions, for which are
    grateful.
    Finally, we thank the participants of the Ramanujan
    Complexes Seminar that took place at the Hebrew University in the winter of 2013.
}

\section*{Notation}

	Throughout, all vector spaces are over $\C$.
	An algebra means a unital $\C$-algebra. Modules
	are assumed to be unital.\footnote{In \cite{Fi15A}, algebras and modules are not
	assumed to be unital.}
	Subalgebras are not required to have the same unity as the ambient algebra. 
	For an algebra $A$, a left $A$-module $V$ and $a\in A$, denote by $a|_V$ the linear operator $[v\mapsto av]\in \End_{\C}(V)$. 
	
	If $V$ is a Hilbert space, then $\sphere{V}$ denotes the unit sphere of $V$.
	If $X$ is a set, then we write $Y\fsubseteq X$ to denote that $Y$ is a finite subset of $X$.
	
\medskip
	
	For a set $X$, we let $\llf(X)$ denote the set of functions $\vphi:X\to \C$ with finite
    support. We endow $\llf(X)$ with the inner product $\Trings{\vphi,\psi}=\sum_{x\in X}\vphi x\cdot\cconj{\psi x}$.
    This makes $\llf(X)$ into a pre-Hilbert space. Its completion is the Hilbert space of square-summable functions
    on $X$, denoted $\ell^2(X)$.
    The vector space $\llf(X)$ admits a standard basis $\{\e_x\}_{x\in X}$
    defined by
    \[
    \e_x(y)=\left\{\begin{array}{ll}
    1 & x=y \\
    0 & x\neq y
    \end{array}\right.\ .
    \]
    If $Y$ is another set and $f:X\to Y$ is any function, then we define $f_*:\llf(X)\to \llf(Y)$
    by $(f_*\vphi)y=\sum_{x\in f^{-1}\{y\}}\vphi(x)$ for all $\vphi\in\llf(Y)$, $y\in Y$.
    In particular, we have
    \[
    f_*\e_x=\e_{f(x)}\qquad\forall\, x\in X\ .
    \]

\medskip
	
	Recall that an $\ell$-group is a locally compact totally disconnected
	Hausdorff topological group. Such groups admit a basis of neighborhoods at
	the identity consisting of compact open subgroups.

\section{$*$-Algebras}
\label{sec:algebras}

	In the basis of our definition of spectrum of simplicial complexes lies
	the spectral theory of \emph{idempotented $*$-algebra}, as developed in \cite[Ch.~2]{Fi15A}.
	We shall restrict here to unital $*$-algebra for the sake of simplicity.
	General $*$-algebras are discussed in \cite{Palmer01}, for instance.
	
\rem{
\medskip

	Throughout, all vector spaces are over $\C$.
	An algebra means a unital $\C$-algebra. Subalgebras are not required to have the same unity as the ambient algebra. Modules
	are assumed to be unital.\footnote{In \cite{Fi15A}, algebras and modules are not
	assumed to be unital.}
	
	For a left $A$-module $V$ and $a\in A$, denote by $a|_V$ the linear operator $[v\mapsto av]\in \End_{\C}(V)$. If $V$ is a Hilbert space, then $\sphere{V}$ denotes the unit sphere of $V$.
}

\subsection{Unitary Representations}
\label{subsec:idempotented-algebras}



	Let $A$ be an algebra.
    An involution on $A$ is a map $*:A\to A$ such that $a^{**}=a$, $(a+b)^*=a^*+b^*$, $(ab)^*=b^*a^*$
    and $(\alpha a)^*=\cconj{\alpha} a^*$ for all $a,b\in A$ and $\alpha\in\C$.
    A \emph{$*$-algebra} is an algebra $A$ equipped with an involution, which is always denoted by $*$.
    
    \begin{example}
    	The commutative algebra $A=\C[X_1,\dots,X_t]$ is a $*$-algebra with respect to the unique involution
    	$*$ satisfying $X_i^*=X_i$ for all $i$.
    \end{example}



    
    A \emph{unitary representation} of $A$ is a Hilbert space $V$ equipped with a left
    $A$-module structure such that
    \begin{enumerate}
        \item[(U1)] $\Trings{au,v}=\Trings{u,a^*v}$ for all $a\in A$ and $u,v\in V$,
        \item[(U2)] for all $a\in A$, the operator $a|_{V}:V\to V$ is  bounded.
    \end{enumerate}
    We say that $V$ is \emph{irreducible}\footnote{
        Some texts use the term \emph{topologically irreducible}.
    } if it is does not have a proper nonzero \emph{closed}  $A$-submodule.
    Let $\Rep[u]{A}$ denote the category whose objects are unitary representations of $A$
    and whose morphisms are continuous $A$-module homomorphisms.
    Morphisms preserving the inner product are called \emph{unitary}.
    We  let $\Irr[u]{A}$ denote the class of irreducible unitary representations of $A$.




\medskip

    Let $\{V_i\}_{i\in I}\subseteq\Rep[u]{A}$. The direct sum $\bigoplus_i V_i$
    admits an obvious inner-product making it into a pre-Hilbert space. The completion of  $\bigoplus_i V_i$
    is denoted $\hat{\bigoplus}_iV_i$.
    If $\sup_{i}\norm{a|_{V_i}}<\infty$ for all $a\in A$, then the diagonal action of $A$ on $\bigoplus_iV_i$
    extends to $\hat{\bigoplus}_iV_i$ and we may regard $\hat{\bigoplus}_iV_i$ as a unitary
    representation of $A$. We denote this by writing $\hat{\bigoplus}_iV_i\in\Rep[u]{A}$.
    When $I$ is finite, we always have $\bigoplus_iV_i=\hat{\bigoplus}_iV_i\in\Rep[u]{A}$.
    
    We write $V_1\leq V_2$ if there is a unitary injective $A$-homomorphism
    from $V_1$ to $V_2$.

\medskip

    We now recall several well-known facts about unitary representations.

    \begin{thm}[Schur's Lemma; {\cite[Th.~2.6]{Fi15A}}]\label{TH:Schur-lemma}
        Let $V\in \Irr[u]{A}$. Then the continuous $A$-endomorphisms of $V$ are $\C\id_V$.
    \end{thm}

    \begin{cor}[{\cite[Cor.~2.7]{Fi15A}}]\label{CR:commutative-algs}
        If $A$ is commutative, then all irreducible unitary representations of $A$
        are $1$-dimensional.
    \end{cor}

    \begin{prp}[{\cite[Prp.~2.8]{Fi15A}}]\label{PR:iso-implies-unitary-iso}
        Let $V,V'\in\Rep[u]{A}$. If there exists a continuous $A$-module isomorphism
        $f:V\to V'$, then there exists a \emph{unitary} $A$-module isomorphism
        $g:V\to V'$. 
    \end{prp}

\subsection{Spectrum and The Unitary Dual}

	Let $A$ be a $*$-algebra. The isomorphism class of $V\in \Rep[u]{A}$ is denoted $[V]$.
	The set of isomorphism classes of irreducible unitary representations of $A$ is called
	the \emph{unitary dual} of $A$ and denoted
	\[
    \udual{A}=\{[V]\where V\in\Irr[u]{A}\}\ .
    \]
    A subset $S$ of $\what{A}$ is called \emph{bounded} if
    $\hat{\bigoplus}_{[V]\in S}V\in\Rep[u]{A}$, or equivalently,
    if $\sup_{[V]\in S}\norm{a|_V}<\infty$ for all $a\in A$.

	\begin{example}
		Consider $A=\C[X_1,\dots,X_t]$ with the unique involution satisfying $X_i^*\cong X_i$.
		By Corollary~\ref{CR:commutative-algs}, irreducible unitary representations of $A$
		are $1$-dimensional. In particular, they are irreducible
		$A$-modules. For $\lambda=(\lambda_1,\dots,\lambda_t)\in\C$,
		let $V_\lambda=A/\Trings{X_1-\lambda_1,\dots,X_t-\lambda_t}$. It
		is well-known that $\{V_\lambda\where \lambda\in\C^t\}$ form a complete
		set irreducible $A$-modules up to isomorphism. 
		However, only those $V_\lambda$ for which $\lambda\in \R^t$
		can be made into unitary representations of $A$. The unitary dual of $A$ is therefore
		in one-to-one correspondence with $\R^n$. Explicitly, the isomorphism $\what{A}\to \R^n$
		is given by $[V]\mapsto (X_1|_V,\dots,X_t|_V)$.
		
		More generally, if $A$ is any commutative $*$-algebra generated as a $*$-algebra
		by $a_1,\dots,a_t$, then the map 
		\[
		[V]\mapsto (a_1|_V,\dots,a_t|_V)\in\C^t
		\]
		gives an embedding of $\udual{A}$ in $\C^t$. See \cite[Rm.~2.23]{Fi15A} for details
		about its image.
	\end{example}
	
	We make $\udual{A}$ into a topological space as follows:
	Let $V\in\Irr[u]{A}$, $v\in \sphere{V}$, $\veps>0$
    and  $F\fsubseteq A$.
    We define
    \[
    N_{V,v,\veps,F}\subseteq \udual{A}
    \]
    to be the set of all isomorphism classes $[U]\in\udual{A}$ for which
    there is $u\in {U}$ such that
    \[
    \Abs{\Trings{av,v}-\Trings{au,u}}<\veps\qquad\forall\, a\in F\ .
    \]
    Note that $u$ is not required to be a unit vector.
    The possible sets $N_{V,v,\veps,F}$ form a subbasis for a topology on $\udual{A}$.
	
\medskip	
	
	Let $[V]\in\udual{A}$ and $V'\in\Rep[u]{A}$. We say that $V$ is weakly contained in $V'$
	and write
	\[
	V\wc V'
	\]
	if the following equivalent conditions hold (cf.\ \cite[Lm.~2.15]{Fi15A}):
	\begin{enumerate}
		\item[(a)] For all $v\in\sphere{V}$, $\veps>0$ and $F\fsubseteq A$,
        there exists $v'\in\sphere{V'}$ such that $\abs{\Trings{av,v}-\Trings{av',v'}}<\veps$
        for all $a\in F$.
        \item[(b)] There exists $v\in\sphere{V}$ such that for all $\veps>0$ and $F\fsubseteq A$,
        there is $v'\in{V'}$ such that $\abs{\Trings{av,v}-\Trings{av',v'}}<\veps$
        for all $a\in F$.
	\end{enumerate}
	For example, if $V\leq V'$, then $V\wc V'$. The converse is false in general.
	The \emph{$A$-spectrum} of $V'$ is a subset of $\udual{A}$ defined as 
	\[
	\Spec_A(V')=\{[V]\in\udual{A}\suchthat V\wc V'\}\ .
	\]
	
	The following proposition shows that when $A$ is commutative and
	generated as a $*$-algebra by $a_1,\dots,a_t$, there is a 
	topological embedding  of $\udual{A}$ in $\C^t$ such that for every
	$V'\in\Rep[u]{A}$, the set $\Spec_A(V')$ corresponds to 
	the common (continuous) spectrum of $(a_1,\dots,a_t)$ on $V$, denoted
	$\Spec(a_1|_{V'},\dots,a_t|_{V'})$. The $A$-spectrum
	is therefore essentially equivalent to the common spectrum
	of $(a_1,\dots,a_t)$.
	
	\begin{prp}[{\cite[Pr.~2.22]{Fi15A}}]\label{PR:unitary-dual-topology}
        Assume $A$ is commutative and generated
        as a $*$-algebra by $a_1,\dots,a_t$. For $V\in \Irr[u]{A}$,
        denote denote by $\lambda_V\in\C^t$ the unique common eigenvalue
        of $(a_1,\dots,a_t)$ on $V$ (cf.\ Corollary~\ref{CR:commutative-algs}).
        Then the map
        \[
        [V]\mapsto \lambda_V: \udual{A}\to\C^t
        \]
        is a topological embedding. In addition, for all $V'\in\Rep[u]{A}$, we have
        \[
        \Spec(a_1|_{V'},\dots,a_t|_{V'})=\{\lambda_V\where [V]\in\Spec_A(V')\}\ .
        \]
        Finally, a subset $S\subseteq\udual{A}$ is bounded if and only if its image in $\C^t$
        is bounded.
    \end{prp}

	The unitary dual of finitely generated 
	non-commutative algebras is not Hausdorff
	in general
	\cite[Ex.~2.24]{Fi15A}.
	
\medskip

	We  mention here several useful results about the $A$-spectrum.
	
	\begin{prp}[{\cite[Pr.~2.28]{Fi15A}}]\label{PR:spectrum-is-closed}
        Let $V'\in\Rep[u]{A}$. Then $\Spec_A(V')$ is bounded and closed in $\udual{A}$.
    \end{prp}
    
    \begin{prp}[{\cite[Cor.~2.30]{Fi15A}}]\label{PR:non-empty-spectrum}
    	Let $0\neq V'\in\Rep[u]{A}$ and let $a\in A$. There
    	is $[V]\in\Spec_A(V')$ such that $\norm{a|_V}=\norm{a|_{V'}}$.
    	In particular, $\Spec_A(V')\neq\emptyset$.
    \end{prp}
	
	\begin{thm}[{\cite[Th.~2.31]{Fi15A}}]\label{TH:direct-sum-spectrum}
        Let $\{V_i\}_{i\in I}$ be a family of unitary representations of $A$ such
        that $\hat{\bigoplus}_iV_i\in\Rep[u]{A}$. Then
        $
        \Spec_A({\hat{\bigoplus}_iV_i})=\overline{\bigcup_i\Spec_A(V_i)}
        $.
    \end{thm}
    
\subsection{Subalgebras}

	Let $A$ be a $*$-algebra. A $*$-subalgebra of $A$ is
	a subalgebra $B$ such that $B^*=B$. If $V\in\Rep[u]{A}$,
	then $BV$ is a unitary
	representation of $B$ (notice that $BV$ is closed
	in $V$ since $B$ has a unity). When the unities of $B$ and $A$
	are the same, we have $BV=V$.
	The following theorem shows that the $A$-spectrum of $V$ determines
	the $B$-spectrum of $BV$. 
	
	\begin{thm}[{\cite[Th.~2.35]{Fi15A}}]\label{TH:subalgebra-spectrum-I}
        For all $V'\in\Rep[u]{A}$, we have 
        \[\Spec_B({BV'})=\{[U]\in\udual{B}\suchthat
        \text{there is $[V]\in\udual{A}$ with $U\leq BV$}\}\ .\]
    \end{thm}
    
    \begin{cor}[{\cite[Cor.~2.37]{Fi15A}}]\label{CR:subalgebra-spectrum-II}
        Let $V'\in\Rep[u]{A}$, and let $a_1,\dots,a_t\in A$
        be elements generating a commutative $*$-subalgebra of $A$.
        Then
        \[\Spec(a_1|_{V'},\dots,a_t|_{V'})=\bigcup_{[V]\in\Spec_A(V')}\Spec(a_1|_{V},\dots,a_t|_{V})\ .
        \]
    \end{cor}

	\begin{thm}[{\cite[\S{}2I]{Fi15A}}]
		Let $e\in A$ be an idempotent with $e^*=e$ and let $f=1-e$.
		Then $eV\in\Irr[u]{eAe}$ for all $V\in\Irr[u]{A}$
		with $eV\neq 0$, and for all $V'\in\Rep[u]{A}$, we have
		\begin{align*}
			& \Spec_{eAe}(eV')=\{[eV]\where [V]\in\Spec_A(V'),\,eV\neq 0\},\\
			& \Spec_{fAf}(fV')=\{[fV]\where [V]\in\Spec_A(V'),\,fV\neq 0\},\\
			& \Spec_A(V')=\{[V]\in\udual{A}\suchthat \text{$[eV]\in \udual{eAe}$ or
			$[fV]\in\udual{fAf}$}\}\ .
		\end{align*}
		In particular,  $\Spec_A(V')$ determines $\Spec_{eAe}(eV')$
		and $\Spec_{fAf}(fV')$, and vice versa.
	\end{thm}
	
	In  general, it is not true that for every $U\in\Irr[u]{eAe}$
	there is $V\in\Irr[u]{A}$ with $U\cong eV$ \cite[Rm.~2.43]{Fi15A}.

\subsection{Pre-Unitary Representations}

	Let $A$ be a $*$-algebra. A \emph{pre-unitary representation} of $A$ is a pre-Hilbert space $V$
    endowed with a left $A$-module structure satisfying conditions (U1) and (U2) of \ref{subsec:idempotented-algebras}.

\medskip

    If $V$ is a pre-unitary representation, then the action of $A$ extends to the completion
    $\quo{V}$, which then becomes a unitary representation. The $A$-spectrum of $V$
    is defined to be the $A$-spectrum of $\quo{V}$.

\section{Simplicial Complexes}
\label{sec:simp-comp}



\subsection{Simplicial Complexes}
\label{subsec:simplicial-complexes}

	Recall that a  \emph{simplicial complex} consists of a non\-empty  set  of finite sets  $X$
    such that subsets of sets in $X$ are also in $X$.
    A partially ordered set $(Y,\leq)$ that is isomorphic to $(X,\subseteq)$ for some simplicial complex $X$ will also be called
    a simplicial complex. 

    The elements of a simplicial complex $X$ are called \emph{cells}.
    It is \emph{locally finite} if every cell in $X$ is contained in finitely many cells.
    We let $X^{(i)}$ denote the sets in $X$ of cardinality
    $i+1$. Elements of $X^{(i)}$ are called $i$-dimensional cells, or just
    $i$-cells.
    The \emph{vertex set} of $X$ is  $\vrt{X}:=\bigcup_{x\in X}x$. By
    abuse of notation, we sometimes refer to elements of $X^{(0)}$ as vertices.
    The dimension of $X$ is the maximal $i$ such that $X^{(i)}\neq \emptyset$.
    
    If $x$ and $x'$ are distinct cells in $X$, then the \emph{combinatorial distance} of $x$ from $x'$,
    denoted $\dist(x,x')$, is the minimal $t\in\N$
    such that there exists a sequence of cells $y_1\dots,y_t$  with
    $x\subseteq y_1$, $x'\subseteq y_t$ and
    $y_i\cap y_{i+1}\neq \emptyset$
    for all $1\leq i<t$. We further set $\dist(x,x)=0$. (This agrees with the combinatorial distance in graphs.)
    The \emph{ball of radius $n$} around $x$,  $\Ball_X(x,n)$, consists of the
    cells in $X$ of distance $n$ or less from $x$. When $X$ is locally finite, all the balls $\Ball_X(x,n)$  ($x\neq\emptyset$)
    are  finite.
    We say that $X$ is \emph{connected} if $X=\bigcup_{n\geq 0} \Ball_X(x,n)$ for some (and hence all) $x\in X- \{\emptyset\}$,
    or equivalently, if $\dist(x,x')<\infty$ for all $x,x'\in X$.

	A morphism of simplicial complexes $f:X\to Y$ consists of a function $f:\vrt{X}\to \vrt{Y}$
    such that for all $i$ and $x\in X^{(i)}$, we have $f(x):=\{f(v)\where v\in x\}\in Y^{(i)}$. The induced
    maps $X^{(i)}\to Y^{(i)}$ and $X\to Y$ are also denoted $f$. 
    A morphism $f:X\to Y$  is a \emph{cover map}
    if it  is a cover map when $X$ and $Y$ are realized as topological spaces
    in the obvious way. This is equivalent to
    saying that $f:\vrt{X}\to \vrt{Y}$ is surjective and  the induced map $f:\{x\in X\suchthat v\in x\}\to \{y\in Y\suchthat
    f(v)\in y\}$ is bijective for all $v\in\vrt{X}$.
    In this case, the \emph{deck transformations} of  $f:X\to Y$ are the
    automorphisms $h$
    of $X$ satisfying $f\circ h=f$.
    We let
    \[
    \catSimp
    \]
    denote the category of locally finite connected simplicial complexes with cover maps as morphisms.
    
    
\medskip

	\emph{Henceforth, simplicial complexes are assumed to be connected and locally
	finite.}
	
\subsection{Orientation}
\label{subsec:orientation}

	Let $X$ be a simplicial complex. An \emph{ordered cell} in $X$ consists of a pair $(x,\leq)$ where $x\in X$ and $\leq$ is a full
    ordering of the vertices of $x$.
    Two orders on $x$ are equivalent if one can be obtained from the other by an even permutation. We denote
    by $[x,\leq]$ the equivalence class of $(x,\leq)$ and call it an \emph{oriented cell}.
    When $|x|>1$, we also write $[x,\leq]^\op$ to denote $x$ endowed with orientation different
    from the one induced by $\leq$. 
    For  $\{v_0,\dots,v_i\}\in X^{(i)}$, let
    \[
    [v_0v_1\dots v_i]=[\{v_0,\dots,v_i\}\, ,\,v_0<v_1<\dots<v_i]\ .
    \]
    The collection of oriented $i$-dimensional
    cells is denoted $\ori{X}^{(i)}$.
    
    For every $i>0$, define
    \begin{align*}
    \Omega_i^\pm(X)&:=\llf(\ori{X^{(i)}})\\
    \Omega_i^-(X)&:=\{\vphi\in\Omega_i^\pm(X)\suchthat \text{$\vphi(\sfx^\op)=-\vphi(\sfx)$
    for all $\sfx\in \ori{X^{(i)}}$}\}\\
    \Omega_i^+(X)&:=\{\vphi\in\Omega_i^\pm(X)\suchthat \text{$\vphi(\sfx^\op)=\vphi(\sfx)$
    for all $\sfx\in \ori{X^{(i)}}$}\}
    \end{align*}
    The inner product on $\llf(\ori{X^{(i)}})$ makes all three spaces
    into pre-Hilbert spaces.
    We further write
    \[
    \Omega_0^\pm(X)=\Omega_0^-(X)=\Omega_0^+(X):=\llf(\ori{X^{(0)}})
    \]
    and endow them with the inner product $\Trings{\vphi,\psi}_{\Omega_0^\pm(X)}=
    2\Trings{\vphi,\psi}_{\llf(\ori{X^{(0)}})}$.
    The space $\Omega_i^-(X)$ (resp.\ $\Omega_i^+(X)$) is the space of $i$-dimensional
    forms (resp.\ anti-forms) on $X$.
    
    Observe that $\Omega_i^\pm$ is a covariant functor from $\catSimp$ to 
    $\catPHil$, the category of pre-Hilbert spaces with (non-continuous) 
    linear maps as morphisms; a morphism $f:X\to Y$ in $\catSimp$
    is mapped to the linear map $f_*:\llf(\ori{X^{(i)}})\to \llf(\ori{Y^{(i)}})$
    determined by $f\e_{\sfx}=\e_{f\sfx}$  for all $\sfx\in\ori{X^{(i)}}$
    ($f\sfx$ is defined in the obvious way; cf.\ the notation section).
    Likewise, $\Omega_i^-$ and $\Omega_i^+$ are subfunctors of $\Omega_i^\pm$.
    Notice that $\Omega_i^+(X)$ is naturally isomorphic to $\llf(X^{(i)})$ (as pre-Hilbert
    spaces). We will therefore occasionally identify $\Omega_i^+(X)$ with $\llf(X^{(i)})$.
    
\medskip

	Recall that the boundary map $\partial_{i+1,X}:\Omega_{i+1}^-(X)\to \Omega_{i}^-(X)$
    and the coboundary map
    $\delta_{i,X}:\Omega_i^-(X)\to\Omega_{i+1}^-(X)$
    are defined by
        \[
        (\partial_{i+1,X}\psi)[v_0\dots v_i]=\sum_{\substack{v\in\vrt{X} \\ \{v,v_0,\dots,v_i\}\in X^{(i+1)}}}\psi[vv_0\dots v_1]
        \]
        \[
        (\delta_{i,X}\vphi)[v_0\dots v_{i+1}]=\sum_{j=0}^{i+1}(-1)^j\vphi[v_0\dots \hat{v}_j\dots v_{i+1}]
        \]
    ($\hat{v}_j$ means omitting the $j$-th entry).
    It is easy to check that $\delta_{i,X}^*=\partial_{i+1,X}$.
    The upper, lower and total $i$-dimensional Laplacians are defined by
        \[
        \Delta_{i,X}^+=\partial_{i+1,X}\delta_{i,X},\qquad\Delta_{i,X}^-=\delta_{i-1,X}\partial_{i,X},\qquad\Delta_{i,X}=\Delta_{i,X}^++\Delta_{i,X}^-\, ,
        \]
    respectively (with the convention that $\Delta_{0,X}^-=0$). 
    It is easy to check that $\partial_{i+1}=\{\partial_{i+1,X}\}_{X\in\catSimp}$
    is a natural transformation from $\Omega_{i+1}^-$ to $\Omega_i^-$,
    and  $\delta_i=\{\delta_{i,X}\}_{X\in\catSimp}$ is a natural
    transformation from $\Omega_{i}^-$ to $\Omega_{i+1}^-$.
    Likewise, $\Delta_i^+$, $\Delta_i^-$ and $\Delta_i$ are natural transformations
    from $\Omega_i^-$ into itself.
    
\subsection{Group Actions}
\label{subsec:quotienst-of-simp-comps}

	Let $G$ be an $\ell$-group (see the notation section).
    By a \emph{$G$-complex} we mean a  (locally finite, connected) simplicial complex $\calX$ on which $G$ acts faithfully via automorphisms
    and such that for all nonempty $x\in\calX$, the stablizer $\Stab_G(x)$ is a compact open subgroup
    of $G$. 
    The $\ell$-groups $G$ for which a given simplicial complex $\calX$ is a $G$-complex
    are characterized in the following proposition.
    
    \begin{prp}[{\cite[Pr.~3.3]{Fi15A}}]\label{PR:characterization-of-G-comps}
        Let $\calX$ be a simplicial complex. Give $\calX$ the discrete topology
        and $\Aut(\calX)$ the topology of pointwise convergence.
        Then:
        \begin{enumerate}
            \item[(i)] $\Aut(\calX)$ is an $\ell$-group and $\calX$ is a $\Aut(\calX)$-complex. Consequently, $\calX$ is a
            $G$-complex for any closed subgroup $G$ of $\Aut(\calX)$.
            \item[(ii)] If $\calX$ is a $G$-complex for an $\ell$-group  $G$, then the action map $G\to \Aut(\calX)$
            is a closed embedding.
        \end{enumerate}
    \end{prp}
    
    A $G$-complex $\calX$ is called \emph{almost transitive} if $G\leftmod\calX$ is finite.
    
    \begin{example}[{\cite[Ex.~3.4]{Fi15A}}]\label{EX:G-complex-building-III}
        Let $\bfG$ be an almost simple algebraic group over a non-archimedean
        local field $F$, let $\bfZ$ be the center of $\bfG$, and let $G=\bfG(F)/\bfZ(F)$.\footnote{
            In general, $\bfG(F)/\bfZ(F)$ is not the same as $(\bfG/\bfZ)(F)$. More precisely, one has
            an exact sequence $1\to \bfZ(F)\to \bfG(F)\to (\bfG/\bfZ)(F)\to \mathrm{H}^1_{\textit{fppf}}(F,\bfZ)$.
        }
        Let $\calB$ be the  \emph{affine Bruhat--Tits
        building} of $\bfG$ (see \cite{Tits79}, \cite{BruhatTits72I}; a more elementary
        treatment in the case $\bfG$ is  classical
        can be found in \cite{AbramNebe02}).
        The building $\calB$ is a  simplicial complex carrying a faithful left $G$-action
        making it into  an almost transitive $G$-complex.
    \end{example}
    
    Let $\calX$ be a $G$-complex
    and
    let $\Gamma$ be a subgroup of $G$.
    The partial order on $\calX$ induces a partial order on $\Gamma\leftmod \calX$ given
    by $\Gamma x\leq \Gamma y$ $\iff$ $\gamma x\subseteq y$ for some $\gamma\in\Gamma$.
    However, $\Gamma\leftmod \calX$ is not a simplicial complex in general
    \cite[Ex.~3.5(ii)]{Fi15A}, and even when this is the case,
    the projection map $x\mapsto\Gamma x:\calX\to \Gamma\leftmod \calX$ may not be a cover map
    \cite[Ex.~3.6]{Fi15A}.
    When both conditions hold, we call $\Gamma\leftmod \calX$  a \emph{$G$-quotient} of $\calX$
    and write 
    \[\Gamma\leq_\calX G\ .\]
    In this case, $\Gamma$ coincides with the group of deck transformations of $\calX\to\Gamma\leftmod \calX$
    \cite[Pr.~3.9]{Fi15A}.

    \begin{prp}[{\cite[Pr.~3.7, Pr.~3.9]{Fi15A}}]
    	Let $\Gamma\leq G$. Then $\Gamma\leftmod\calX$ is a simplicial complex
    	if and only if
    	\begin{enumerate}
    		\item[(C1)] $\{\Gamma u_1,\dots,\Gamma u_t\}=\{\Gamma v_1,\dots,\Gamma v_s\}$
            implies $\Gamma\{u_1,\dots,u_t\}=\Gamma\{v_1,\dots,v_s\}$ for all $\{u_1,\dots,u_t\},\{v_1,\dots,v_s\}\in \vrt{\calX}$. 
        \end{enumerate}
        In this case, $\dim (\Gamma\leftmod \calX)=\dim \calX$ and the map $x\mapsto \Gamma x: \calX\to \Gamma\leftmod \calX$ is a morphism of
        simplicial
        complexes. The latter map is a cover map if and only if
        \begin{enumerate}
            \item[(C2)] $\Gamma\cap\Stab_G(v)=\{1_G\}$ for all $v\in\vrt{\calX}$.
        \end{enumerate}
        In this case, $\Gamma$ acts freely on $\calX-\{\emptyset\}$ and
        $\Gamma$ is discrete in $G$.
    \end{prp}
    
    The proposition can be used to prove the following elegant criterion for
    when $\Gamma\leq_{\calX}G$.
    
    \begin{cor}[{\cite[Cor.~3.10]{Fi15A}}]
    	$\Gamma\leq_{\calX} G$ $\iff$ $\dist(v,\gamma v)>2$ for all $1\neq\gamma\in \Gamma$ and $v\in \vrt{\calX}$.
    \end{cor}
    
    \begin{cor}[{\cite[Cor.~3.11]{Fi15A}}]\label{CR:Gamma-is-discrete}
        If $\Gamma\leq_\calX G$, then $\Gamma'\leq_\calX G$ for all $\Gamma'\leq\Gamma$,
        and $g^{-1}\Gamma g\leq_{\calX}G$ for all $g\in G$.
    \end{cor}
    
    We finish  with noting that if $\Gamma\leq_{\calX}G$,
    then  $\Gamma\leftmod\calX$ is finite if and only if 
    $\Gamma$ is cocompact in $G$, and in this case $G$ is unimodular
    {\cite[Pr.~3.13]{Fi15A}}.

\section{Spectrum in Simplicial Complexes}
\label{sec:spectrum}

	In this section, we introduce our notion of spectrum of simplicial complexes.
	Our definition is  general and we shall demonstrate how it specializes
	to other notions of spectrum considered in the literature, e.g.\ the spectrum
	of $k$-regular graphs.
	The idea is to choose  a  $*$-algebra of operators whose spectrum we wish to investigate
	and apply the spectral theory of Section~\ref{sec:algebras}. 
	
\subsection{Associated Operators}

	Let $\catC$ be a subcategory of $\catSimp$, the category of locally
	finite connected simplicial complexes with cover maps as morphisms (\ref{subsec:simplicial-complexes}),
	and let $F$ be a (covariant) functor from $\catC$ to $\catPHil$, the category
	of pre-Hilbert spaces with (non-continuous) linear maps as morphisms.

	\begin{dfn}[{\cite[Def.~4.5, \S{}4C]{Fi15A}}]
		An \emph{associated operator} of $(\catC,F)$, or just a \emph{$(\catC,F)$-operator}, 
		is a natural transformation $a=\{a_X\}_{X\in\catC}:F\to F$
        that admits a dual, i.e.\ a natural transformation $a^*=\{a^*_X\}_{X\in\catC}:F\to F$
        such that for all $X\in\catC$, the operator $a^*_X$ is the dual of $a_X:FX\to FX$ relative to the inner product
        of $FX$.
        
        The collection of all $(\catC,F)$-operators is denoted
        \[
        A(\catC,F)\ .
        \]
    	An \emph{algebra of $(\catC,F)$-operators} is a subset of $A(\catC,F)$
    	that is closed under addition, composition, scaling by elements of $\C$, and
    	taking the dual transformation.
	\end{dfn}
	
	The definition still works if one replaces $\catSimp$ with the category
	of locally
	finite connected simplicial complexes with \emph{arbitrary} morphisms.

	\begin{example}[{\cite[Ex.~4.3, \S{}4B]{Fi15A}}]
		Take $\catC=\catSimp$. We observed in \ref{subsec:orientation}
		that $\Omega_i^-:\catSimp\to \catPHil$ is a functor
		and that the $i$-dimensional Laplactian $\Delta_i:\Omega_i^-\to \Omega_i^-$
		is a self-dual natural transformation. Thus, $\Delta_i$
		is a $(\catSimp,\Omega_i^-)$-operator,
		and $\C[\Delta_i]$, the $\C$-algebra  spanned by $\Delta_i$ in 
		$A(\catSimp,\Omega_i^-)$, is an algebra of $(\catSimp,\Omega_i^-)$-operators.
	\end{example}
	
	In the next examples, we identify $\Omega_i^+(X)$ with $\llf(X^{(i)})$
	(cf.~\ref{subsec:orientation}).
	
	\begin{example}[{\cite[Ex.~4.1, \S{}4B]{Fi15A}}]\label{EX:adjacency-operator-tree}
		Let $\catC\subseteq\catSimp$ be the full subcategory of $k$-regular graphs.
		For $X\in\catC$, the \emph{vertex adjacency operator} 
		$a_{0,X}:\Omega_0^+(X)\to\Omega_0^+(X)$
		is defined by
        \[(a_{0,X}\vphi)u=\sum_{v}\vphi u\qquad\forall\, \vphi\in\Omega^+_0(X),\, u\in X^{(0)}\,,\]
        where the sum is taken over all $v\in X^{(0)}$ connected by an edge to $u$.
		Then $a_0:=\{a_{0,X}\}_{X\in\catC}$ is an associated operator of $(\catC,\Omega_0^+)$
		and $\C[a_0]$ is an algebra of $(\catC,\Omega_0^+)$-operators. In fact,
		$\C[a_0]=A(\catC,\Omega_0^+)$ \cite[Ex.~4.20]{Fi15A} (this
		is false for a general subcategory $\catC$).
	\end{example}
	
	\begin{example}[{\cite[Ex.~4.2]{Fi15A}}]\label{EX:higher-dimensional-adj-ops}
		Let $X$ be a simplicial complex and assume $i,j$ satisfy $0\leq i<j\leq 2i+1$.
        Define $a_{i;j,X}:\Omega^+_i(X)\to \Omega^+_i(X)$ by
        \[(a_{i;j,X}\psi)x=\sum_{\substack{y\in X^{(i)}\\ x \cup y\in X^{(j)}}}\vphi y\qquad\forall\, \vphi\in\Omega^+_i(X),\, x\in X^{(i)}\ .\]
        That is, the evaluation of $a_{i;j,X}\vphi$  at an $i$-cell $x$ adds the values of $\vphi$ on $i$-cells $y$ whose union with $x$ is a $j$-cell.
        (In the notation of Example~\ref{EX:adjacency-operator-tree}, we have $a_{0,X}=a_{0;1,X}$.)
        Then $a_{i;j}$ is a $(\catSimp,\Omega_i^+)$-operator and the algebra
        spanned by $a_{i;i+1},\dots,a_{i;2i+1}$ is an algebra of $(\catSimp,\Omega_i^+)$-operators.
        It is not commutative when $i>0$.
	\end{example}
	
\subsection{Spectrum}
\label{subsec:spectrum}
	
	Let $A$ be algebra of $(\catC,F)$-operators and let $X\in\catC$.
	Then $A$ acts on $FX$ via $a\cdot v=a_Xv$ ($a=\{a_X\}_{X\in\catC}\in A$,
	$v\in FX$; the action is not necessarily unital).
	If the elements of $A$ act continuously on $FX$,
	then $AFX$ is a pre-unitary representation of $A$,
	and
	we define the \emph{$A$-spectrum} of $X$ to be
	\[
	\Spec_A(X)=\Spec_A(AFX)
	\]
	(if $A$ contains $\id_F:F\to F$, we have $AFX=FX$). When
	$A=A(\catC,F)$ we also write $\Spec_{\catC,F}(X)$ or $\Spec_F(X)$
	instead of $\Spec_A(X)$.
	
	When $\catC$ is understood from the context, we call
	$\Spec_{\Omega_i^+}(X)$ (resp.\ $\Spec_{\Omega_i^-}(X)$, $\Spec_{\Omega_i^\pm}(X)$)
	the non-oriented (resp.\ oriented, full) \emph{$i$-dimensional spectrum},
	and denote it by $\Spec_i(X)$ (resp.\ $\Spec_{-i}(X)$, $\Spec_{\pm i}(X)$).
	There is no ambiguity for $i=0$ since $\Omega_0^+=\Omega_0^-=\Omega_0^\pm$.
	
	\begin{example}[{\cite[Ex.~4.8]{Fi15A}}]\label{EX:k-regular-graphs-basic-spectrum}
        Let $\catC$ be the category of $k$-regular graphs,
        and write $A_0=\Alg{\catC,\Omega^+_0}{}$. It can be checked    that
        \[
        A_0=\C[a_0]
        \]
        where $a_0$ is the vertex adjacency operator of Example~\ref{EX:adjacency-operator-tree}.
        By Proposition~\ref{PR:unitary-dual-topology}, we can identify $\udual{A}_0$
        with a subset of $\C$ (this set is $\R$, in fact) such
        that for all $V\in\Rep[u]{A_0}$, the set $\Spec_{A_0}(V)$
        corresponds to $\Spec(a_0|_V)$. Therefore, for all $X\in\catC$, the 
        datum of $\Spec_{A_0}(X)$
        is equivalent to the spectrum of the vertex adjacency operator of $X$.
        Otherwise stated,
        the (non-oriented) $0$-dimensional spectrum is essentially the same as the 
        usual spectrum of $k$-regular graphs.
        
        Likewise, the non-oriented $1$-dimensional spectrum 
        turns out to be equivalent
        to the spectrum of the edge adjacency operator (for our particular choice
        of $\catC$).
    \end{example}
	
	\begin{example}\label{EX:zero-dim-spec-of-Bd}
		Let $F$ be a local non-archimedean field, let $G=\nPGL{F}{d}$
		and let $\calB_d(F)$ be the affine Bruhat-Tits building of $G$ (see
		\cite[\S{}2]{LubSamVi05} or \ref{subsec:building} below). Let $\catC$ be the collection $\{\Gamma\leftmod\calB_d(F)
		\where \Gamma\leq_{\calB_d(F)} G\}$. Then $\catC$ can be made into a category
		(cf.\ Definition~\ref{DF:G-complex-category} below) such that the (non-oriented) $0$-dimensional
		spectrum of complexes in $\catC$ is  equivalent to the spectrum of quotients
		of $\calB_d(F)$ as defined by Lubotzky, Samuels and Vishne in \cite{LubSamVi05}.
		We give more details about this in \ref{subsec:building}.
    \end{example}
    
    If $A$ is an algebra of $(\catC,F)$-operator containing an algebra
    $B$ of $(\catC,F)$-operators, then the $A$-spectrum (when defined) determines
    the $B$-spectrum, thanks to Theorem~\ref{TH:subalgebra-spectrum-I}.
    In addition,  by Corollary~\ref{CR:subalgebra-spectrum-II}, for all $a_1,\dots,a_t\in A$ such that $a_1,\dots,a_t$, $a_1^*,\dots,a_t^*$
    commute, the $A$-spectrum of $X\in\catC$ (when defined) determines
    the common spectrum of $(a_{1,X},\dots,a_{t,X})$.
    
\subsection{Elementary Functors}

	We now restrict our attention to special families of subcategories $\catC\subseteq\catSimp$
	and functors $F:\catC\to \catPHil$ arising from almost transitive $G$-complexes.
	Henceforth, $G$ is an $\ell$-group and $\calX$ is an almost transitive $G$-complex.
	
	\begin{dfn}[{\cite[Def.~4.6]{Fi15A}}]\label{DF:G-complex-category}
        We define the subcategory
        \[
        \catC=\catC(G,\calX)\subseteq\catSimp
        \]
        as follows: The objects of $\catC$ are $\{\Gamma\leftmod\calX\where \Gamma\leq_{\calX} G\}$
        (see~\ref{subsec:quotienst-of-simp-comps}), where $1\leftmod\calX$ is identified with $\calX$.
        The morphisms of $\catC$ are given as follows:
        \begin{itemize}
            \item For all $\Gamma\leq_{\calX}G$, set $\Hom_{\catC}(\calX,\Gamma\leftmod\calX)=\{p_\Gamma\circ g\where
            g\in G\}$, where $p_\Gamma$ is the quotient map $x\mapsto \Gamma x:\calX\to \Gamma\leftmod\calX$.
            \item For all $1\neq\Gamma'\leq\Gamma\leq_{\calX} G$, set $\Hom_{\catC}(\Gamma'\leftmod\calX,\Gamma\leftmod\calX)=\{p_{\Gamma',\Gamma}\}$
            where
            $p_{\Gamma',\Gamma}$ is the quotient map $\Gamma'x\mapsto \Gamma x:\Gamma'\leftmod\calX\to \Gamma\leftmod\calX$.
            \item All other $\Hom$-sets are empty.
        \end{itemize}
        In particular, $\End_{\catC}(\calX)=G$.
    \end{dfn}
    
    We call $\catC(G,\calX)$ the category of $G$-quotients of $\calX$.
    
    \begin{dfn}\label{DF:elementary-functor}
		A functor $F:\catC(G,\calX)\to \catPHil$ is \emph{elementary}
		if there exists a covariant functor $S:\catC(G,\calX)\to\catSet$ such
		that
		\begin{enumerate}
			\item[{\rm(E1)}] There is a unitary natural
			isomorphism $\llf\circ S\cong F$.
			\item[{\rm(E2)}] For all $x\in S\calX$, the group
			$\Stab_G(x)$ is compact open in $G$ and
			contained in the stabilizer of a nonempty cell
			in $\calX$ (the action of $G$ on $S\calX$ is via $S$).
			\item[{\rm(E3)}] For all $\Gamma\leq_{\calX} G$, the map $\Gamma\leftmod S\calX\to S(\Gamma\leftmod\calX)$ given
            $\Gamma x\mapsto (Sp_\Gamma)x$
			(notation as in Definition~\ref{DF:G-complex-category}) is an isomorphism.
			\item[{\rm(E4)}] $G\leftmod S\calX$ is finite.
		\end{enumerate}
		The functor $F$ is called \emph{semi-elementary} if there is another
		functor $F':\catC(\calX,G)\to \catPHil$ such that $F\oplus F'$ is
		elementary.
	\end{dfn}
	
	\begin{example}[{\cite[Ex.~4.12]{Fi15A}}]
		The functors $\Omega_i^\pm$ and $\Omega_i^+$  are elementary.
		Indeed,
		take $S$ to be $X\mapsto \ori{X^{(i)}}$ and $X\mapsto X^{(i)}$, respectively
		(to  verify (E3), use \cite[Pr.~3.16]{Fi15A}).
		The functor  $\Omega_0^-$ is also elementary.

		For $i>0$, the
		functor $\Omega_i^-$ is semi-elementary since
		$\Omega_i^\pm=\Omega_i^+\oplus\Omega_i^-$.
	\end{example}
	
	For semi-elementary $F:\catC(G,\calX)\to \catPHil$, it is possible to
	give an alternative description of $A(\catC(G,\calX),F)$.
	In addition, the $F$-spectrum is always defined, i.e.\ $A(\catC(G,\calX),F)$
	acts continuously on $FX$ for all $X\in\catC(G,\calX)$. 
	
	Notice that 
	$F\calX$ admits an obvious left $G$-action since $\End_{\catC(G,\calX)}(\calX)=G$.
	
	\begin{thm}[{\cite[Th.~4.15]{Fi15A}}]\label{TH:Adj-algebra-descrition}
        Let  $\catC=\catC(G,\calX)$, let $F:
        \catC\to \catPHil$ be semi-elementary,
        and write $A=\Alg{\catC,F}{}$. Then:
        \begin{enumerate}
            \item[(i)] The map $\{a_X\}_{X\in\catC}\mapsto a_{\calX}:A\to \End_G(F\calX)$ is an isomorphism of 
            $*$-algebras, where the involution on $\End_G(F\calX)$ is given by taking the dual with respect to the inner product on $F\calX$.
            \item[(ii)] For every $a\in A$, there is $M=M(a)\in\R_{\geq0}$ such that $\norm{a|_{FX}}\leq M$
            for all $X\in \catC$. In particular, $\Spec_A(X)$ is defined for all $X\in\catC$.
        \end{enumerate}
    \end{thm}
    
    Let $\catC=\catC(G,\calX)$, let $F:\catC\to\catPHil$ be elementary,
    and write $F=\llf \circ S$ as in Definition~\ref{DF:elementary-functor}.
    Theorem~\ref{TH:Adj-algebra-descrition} allows us to identify
    $\End_G(F\calX)$ with $A:=\Alg{\catC,F}{}$. For $a\in\End_G(F\calX)$
    and $\Gamma\leftmod \calX\in\catC$, the action of $a$ on $F(\Gamma\leftmod\calX)$
    can be described as follows: Let $x\in S\calX$
    and write $a\e_x=\sum_{y\in S\calX}\alpha_y \e_y$ with $\{\alpha_y\}_y\subseteq \C$
    (see the notation section). Then, upon identifying $S(\Gamma\leftmod \calX)$
    with $\Gamma\leftmod S\calX$ as in (E3), we have
    \[
    a\e_{\Gamma x}=\sum_{y\in S\calX}\alpha_y \e_{\Gamma y}
    \]
	(see the proof of {\cite[Th.~4.15]{Fi15A}}).
	
\medskip

	Some examples demonstrating how
	Theorem~\ref{TH:Adj-algebra-descrition} can be used to determine 
	$\Alg{\catC(G,\calX),F}{}$ can be found in \cite[\S{}4D]{Fi15A}.
	It also yields the following result.

	\begin{prp}\label{PR:subgroup-quotient}
    	Let $\Gamma'\leq\Gamma\leq_{\calX}G$ with $[\Gamma:\Gamma']<\infty$,
    	let $F:\catC(G,\calX)\to\catPHil$ be semi-elementary, and let $A$ be an algebra
    	of $(\catC(G,\calX),F)$-operators.
    	Then $\Spec_A(\Gamma\leftmod\calX)\subseteq\Spec_A(\Gamma'\leftmod\calX)$.
    \end{prp}
    
\subsection{Further Remarks}

	Let $\catC$ be a subcategory of $\catSimp$
	and let $F,F':\catC\to \catPHil$ be functors.
	It sometimes happen that the there are dependencies between
	the $F$-spectrum and the $F'$-spectrum. For example, the full $i$-dimensional spectrum
	determines the non-oriented and oriented $i$-dimensional spectra, and vice versa.
	See \cite[\S{}4E, \S{}4F]{Fi15A} for an extensive discussion about this.
	
	We also note that if $X,X'\in\catC$ are isomorphic
    simplicial complexes, then $\Spec_{A}(X)$ and $\Spec_{A}(X')$ may still differ,
    because $X$ and $X'$ may not be isomorphic in $\catC$. 
	In some cases, this problem is only ostensible \cite[Pr.~4.32]{Fi15A}. However,
	this issue does occur for the spectrum of $\nPGL{F}{d}$-quotients of the affine Bruhat--Tits
	building of $\nPGL{F}{d}$ as defined in \cite{LubPhiSar88} (cf.\ Example~\ref{EX:zero-dim-spec-of-Bd});
	see \cite[Rm.~4.33]{Fi15A} for more details.   
    
\section{Optimal Spectrum}
\label{sec:opt-spec}

	For this section, let $\calX$ be an almost transitive $G$-complex, let
	$\catC=\catC(G,\calX)$,  let $F:\catC\to \catPHil$
	be semi-elementary, and let $A$ be an algebra of $(\catC,F)$-operators.
	A theorem in the spirit of the Alon--Boppana Theorem applies to 
	the $A$-spectrum of $G$-quotients of $\calX$.
	
	\begin{thm}[{\cite[Th.~5.1]{Fi15A}}]\label{TH:Alon-Boppana-I}
        Let  $\{\Gamma_\alpha\}_{\alpha\in I}$ be a family
        of subgroups of $G$ such that $\Gamma_\alpha\leq_{\calX} G$ for all $\alpha\in I$,
        and one of the following
        conditions, which are equivalent, is satisfied:
        \begin{enumerate}
        	\item[(1)] For every compact    $C\subseteq G$ with $1\in C$, there exist $\alpha\in I$
        	and $g\in G$ such that
        	such that $C\cap g^{-1}\Gamma_\alpha g=1$
        	\item[(2)] For every $n\in\N$,
        	there exists $\alpha\in I$ and $v\in\vrt{\calX}$
        	such that quotient map $\calX\to\Gamma_\alpha\leftmod\calX$
        	is injective on the  ball $\Ball_\calX(v,n)$.
        \end{enumerate}
        Then
        \[
        \overline{\bigcup_{\alpha}\Spec_{A}(\Gamma_\alpha\leftmod \calX)}\supseteq\Spec_{A}(\calX)\ .
        \]
    \end{thm}
    
    \begin{example}[{\cite[Ex.~5.3]{Fi15A}}]
        Let be $\Gamma_1\supseteq \Gamma_2\supseteq\Gamma_3\supseteq\dots$ be subgroups of $G$ satisfying
        $\bigcap_n\Gamma_n=1$ and $\Gamma_1\leq_{\calX} G$. Then Theorem~\ref{TH:Alon-Boppana-I} applies
        to $\{\Gamma_n\}_{n\in\N}$.
   	\end{example}
	
	Together with Proposition~\ref{PR:unitary-dual-topology} 
	the theorem implies:
	
	\begin{cor}[{\cite[Cor.~5.4]{Fi15A}}]\label{CR:Alon-Boppana-comm}
        Let  $\{\Gamma_\alpha\}_{\alpha\in I}$ be as
        in Theorem~\ref{TH:Alon-Boppana-I},
        and suppose $a_1,\dots,a_t\in\Alg{\catC,F}{}$ satisfy $a_ia_j=a_ja_i$ and
        $a_ia_j^*=a_j^*a_i$ for all $i,j$.
        Write $V_\alpha=F(\Gamma_\alpha\leftmod \calX)$ and $V=F\calX$. Then
        \[
        \overline{\bigcup_{\alpha}\Spec(a_1|_{V_\alpha},\dots,a_t|_{V_\alpha})}\supseteq\Spec(a_1|_{V},\dots,a_t|_{V})\ .
        \]
    \end{cor}	
	
	When $G=\nPGL{F}{d}$ for a local non-archimedean field $F$, $\calX$
	is the affine Bruhat--Tits building of $G$,  $F=\Omega_0^+$, and 
	$a_1,\dots,a_{d-1}$ are \emph{Hecke operators} (the definitions are recalled in \ref{subsec:building} below), Corollary~\ref{CR:Alon-Boppana-comm} 
	is a result of  Li \cite[Th.~4.1]{Li04}.
	
\medskip

	We proceed with describing the \emph{trivial $A$-spectrum}.
	These are special points in the unitary dual $\udual{A}$ that occur
	in the $A$-spectrum of finite $G$-quotients of $\calX$.
	
	Let $N$ be an open finite index subgroup of $G$. Then $F\calX$
	is a left $N$-module and we can form the space of $N$-coinvariants
	\[(F\calX)_N=F\calX/\Span\{v-gv\where v\in F\calX,\, g\in N\}\ .\]
	Since $A$ acts on $F\calX$ via $N$-equivariant linear maps,
	the action of $A$ descends to $(F\calX)_N$. It turns out that
	$A(F\calX)_N$ can be endowed with an inner product making it into
	a unitary representation of $A$, and the unitary isomorphism class
	of $A(F\calX)_N$ is independent of the inner product \cite[Lm.~5.6]{Fi15A}.
	We write
	\[
	\frakT_{A,N}=\Spec_A(A(F\calX)_N)
	\]
	and
	\[
	\frakT_A=\bigcup_N\frakT_{A,N}
	\]
	where $N$ ranges over the open finite-index subgroups of $G$.
	
	The set $\frakT_A$ is called the \emph{trivial $A$-spectrum}. The name is justified
	by the following two results, stating that the question of which points of $\frakT_A$
	occur in $\Spec_A(\Gamma\leftmod \calX)$ is often a matter of which open finite index
	subgroups of $G$ contain $\Gamma$.
	
	\begin{prp}[{\cite[Pr.~5.7]{Fi15A}}]\label{PR:trivial-spectrum-containment}
        Let $\Gamma\leftmod\calX$ be a finite $G$-quotient
        of $\calX$ such that $\Gamma\leq N$. Then $\frakT_{A,N}\subseteq \Spec_A(\Gamma\leftmod\calX)$.
        Furthermore, if $N'\leq N$ is open of finite index, then $\frakT_{A,N}\subseteq \frakT_{A,N'}$.
    \end{prp}
    
    \begin{prp}[{\cite[Pr.~6.34]{Fi15A}}]\label{PR:converse-trivial-spectrum-containment}
    	Assume $F$ is elementary, $A=\Alg{\catC,F}{}$, and write $F=\llf\circ S$
    	as in Definition~\ref{DF:elementary-functor}.
    	Let $x_1,\dots,x_t$ be representatives for the $G$-orbits
    	in $S\calX$ and write $K_n=\Stab_G(x_n)$ ($1\leq n\leq t$). 
		Let $\Gamma\leq_{\calX}G$ be such that $\Gamma\leftmod\calX$ is finite,
		and let $N$ be a \emph{normal} finite-index subgroup of $G$ such that $G/N$ is
		abelian. Then:
        \begin{enumerate}
        	\item[(i)] $\frakT_{A,N}=\bigcup_{n=1}^t\frakT_{A,NK_n}$.
        	\item[(ii)] If $N$ contains $K_n$ for some $n$,
        	then $\Gamma\subseteq N$ $\iff$ $\frakT_{A,N}\subseteq \Spec_A(\Gamma\leftmod\calX)$.
        \end{enumerate} 
    \end{prp}
    
    We note that in many important cases, $G$ contains a minimal open finite index subgroup
    $N$ such that $G/N$ is abelian, in which case Proposition~\ref{PR:converse-trivial-spectrum-containment} 
    applies. For example, this holds when $\calX$ is a $k$-regular tree
    and $G=\Aut(\calX)$ (\cite{Tits70} or \cite{Moller91}), or when $\calX$ is the affine Bruhat--Tits building of
    an almost simple  algebraic group $\bfG$ over a local field $F$ and $G=\im(\bfG(F)\to\Aut(\calX))$.
    
\medskip

	Recall from the introduction that it is of interest to construct
	infinite families of finite $G$-quotients of $\calX$ such that
	their $A$-spectrum is ``as small as possible'', since this is likely
	to manifest in good combinatorial properties. The previous discussion
	suggests the following definition.
	
	\begin{dfn}[{\cite[Def.~5.14]{Fi15A}}]\label{DF:Ramanujan-quotient}
        A $G$-quotient $\Gamma\leftmod\calX$ is \emph{$A$-Ramamnujan} if
        \[
        \Spec_A(\Gamma\leftmod\calX)\subseteq \frakT_A\cup \Spec_A(\calX)\ .
        \]
    \end{dfn}
    
    The
    $A$-Ramanujan $G$-quotients of $\calX$ can be regarded as those quotients whose spectrum
    is as small as one might expect of a decent infinite
    family of $G$-quotients. Alternatively, when finite, they can be regarded as spectral approximations
    of $\calX$.
    
    Here are several possible specializations of Definition~\ref{DF:Ramanujan-quotient}:
    \begin{itemize}
        \item $\Gamma\leftmod\calX$ is \emph{$F$-Ramanujan} if it is $\Alg{F,\catC}{}$-Ramanujan.
        \item $\Gamma\leftmod\calX$ is \emph{Ramanujan in dimension $i$} if it is $\Omega_i^+$-Ramanujan.
        \item $\Gamma\leftmod\calX$ is \emph{completely Ramanujan} if it is $F$-Ramanujan for any
        semi-elementary functor $F:\catC\to \catPHil$.
        (By Proposition~\ref{PR:Ramanujan-well-behaved} below, it is enough
        to check this when $F$ is elementary.)
    \end{itemize}

    \begin{example}[{\cite[Ex.~5.15]{Fi15A}}]\label{EX:classical-Ramamanujan}
        (i) Take $\calX$ to be a $k$-regular tree, $G=\Aut(\calX)$,  $F=\Omega_0^+$ and 
        $A=A(\catC,F)$.
        As explained in Example~\ref{EX:k-regular-graphs-basic-spectrum},
        the $F$-spectrum can be canonically identified with the spectrum of the vertex adjacency operator.
        Careful analysis (see \cite[Ex.~5.12]{Fi15A}) shows that under
        this identification $\frakT_A=\{\pm k\}$ --- same as the usual trivial spectrum of
        $k$-regular graphs --- and that $\Spec_A(\calX)=[-2\sqrt{k-1},2\sqrt{k-1}]$
        (\cite[p.~252, Apx.~3]{Sunada88},
        for instance).
        Thus, a $k$-regular graph is Ramanujan in dimension $0$ (or $\Omega_0^+$-Ramanujan)
        as a $G$-quotient of $\calX$
        if and only if it is Ramanujan in the classical sense.
        In fact, we will  see later
        that $k$-regular graphs  are Ramanujan in dimension $0$ if and only if they are completely Ramanujan.

        (ii) Let $F$ be a non-archimedean local field, let $G=\nPGL{F}{d}$ and let $\calX$
        be the affine Bruhat--Tits building of $G$. It can be checked that a $G$-quotient
        of $\calX$ is Ramanujan in dimension $0$ if and only if it is Ramanujan in the sense of Lubotzky,
        Samuels and Vishne \cite{LubSamVi05}; see \cite[Ex.~5.13, Ex.~5.15(ii)]{Fi15A}.
		%
    \end{example}
    
    The $A$-Ramanujan property behaves well as $A$
    and $\Gamma$ vary.

    \begin{prp}[{\cite[Pr.~5.16]{Fi15A}}]\label{PR:Ramanujan-well-behaved}
        Let $F,F':\catC\to\catPHil$ be semi-elementary functors,  let $A$ be an algebra of $(\catC,F)$-operators,
        and let  $\Gamma'\leq\Gamma\leq_{\calX}G$. Then:
        \begin{enumerate}
            \item[(i)] $\Gamma\leftmod\calX$ is $F\oplus F'$-Ramanujan if and only if
            $\Gamma\leftmod\calX$ is $F$-Ramanujan and $F'$-Ramanujan.
            \item[(ii)] If $\Gamma\leftmod\calX$ is $A$-Ramanuajan, then $\Gamma\leftmod\calX$
            is $B$-Ramanujan for any  $*$-subalgebra $B\subseteq A$.
            \item[(iii)] If $[\Gamma:\Gamma']<\infty$ and $\Gamma'\leftmod\calX$
            is $A$-Ramanujan, then $\Gamma\leftmod\calX$ is $A$-Ramanujan.
        \end{enumerate}
    \end{prp}
    
    Existence of $A$-Ramanujan $G$-quotients  cannot be guaranteed in general. 
    This follows implicitly from  \cite{LubNagn98}.
    Rather, it is more  reasonable to hope that under certain
    assumptions, finite $G$-quotients of $\calX$ would admit \emph{$A$-Ramanujan covers};
    see \cite[Rm.~5.17]{Fi15A} for further discussion.
    
\section{Representation Theory}
\label{sec:rep-thy}

	Unless otherwise indiciated,
	$G$ is a unimodular $\ell$-group, $\mu_G$ is a fixed
	Haar measure of $G$,  $\calX$ is an almost transitive $G$-complex,
	and
    $F:\catC(G,\calX)\to \catPHil$ is an elementary functor.    
    We shall give a criterion for when $\Gamma\leftmod\calX$ is $F$-Ramanujan
    which is phrased in terms of certain unitary representation of $G$.
    
\subsection{Unitary Representations}
\label{subsec:unitary-reps-G}

	As usual, a \emph{unitary representation} of $G$ is a Hilbert space
    $V$ carrying a $G$-module structure such that the action
    $G\times V\to V$ is continuous and
    $\Trings{gu,gv}=\Trings{u,v}$ for all $u,v\in V$, $g\in G$.
    The representation $V$ is \emph{irreducible} if $V$ does not contain closed
    $G$-submodules other than $0$ and $V$. 
    The category of unitary representations of $G$ with continuous $G$-homomorphisms
    is denoted $\Rep[u]{G}$,
    and the class of irreducible representations is denoted $\Irr[u]{G}$. The
    \emph{unitary dual} of $G$, denoted $\what{G}$, is the collection of unitary isomorphism
    classes in $\Irr[u]{G}$.
    
    For every $K\leq G$, we write
    $
    V^K=\{v\in V\suchthat kv=v~\text{for all}~ k\in K\}
    $.
    
    \begin{example}[{\cite[Ex.~6.1, Ex.~6.2]{Fi15A}}]\label{EX:right-reg-rep}
    	Let  $\Gamma$ be a discrete subgroup of $G$. Since $G$ is unimodular, the left coset space 
    	$\Gamma\leftmod G$
    	admits a right $G$-invariant measure $\mu_{\Gamma\leftmod G}$, unique up to scaling.
    	Then $\LL{\Gamma\leftmod G}$ is a unitary representation of $G$
    	with respect to the \emph{left} $G$-action
    	given by 
    	\[
        (g\vphi)x=\vphi(xg)\qquad\forall\, g,x\in G,\, \vphi\in\LL{\Gamma\leftmod G}\ .
    	\]
		When $\Gamma=1$, we have $\mu_{\Gamma\leftmod G}=\mu_G$ (up to scaling),
		and $\LL{\Gamma\leftmod G}=\LL{1\leftmod G}$ is the \emph{right regular representation} of $G$.
	\end{example}

	Let $V\in \Irr[u]{G}$ and $V'\in\Rep[u]{G}$.
    Recall that $V$ is \emph{weakly contained} in $V'$, denoted $V\wc V'$, if
    for all $v\in \sphere{V}$, $\veps>0$ and compact $C\subseteq G$, there exists $v'\in \sphere{V'}$
    such that
    \[
    \abs{\Trings{gv,v}-\Trings{gv',v'}}<\veps\qquad\forall g\in C\ .
    \]
    We further write
    \[
    \Spec_G(V')=\{[V]\in\udual{G}\suchthat V\wc V'\}\ .
    \]

    The representation $V$ is called \emph{tempered} if it weakly contained in $\LL{1\leftmod G}$, the right regular
    representation of $G$.
    (See \cite[\S2.4]{Oh02} for  equivalent definitions of temperedness when $G$ is the group of
    points of a reductive algebraic group over a local field.)
    
\medskip

	Finally, an irreducible representation $V\in\Irr[u]{G}$ is said to have
	\emph{finite action} if one of the following  conditions,
	which are equivalent \cite[Lm.~6.19]{Fi15A}, hold:
	\begin{enumerate}
            \item[(a)] The image of $G$ in the unitary group of $V$ is finite.
            \item[(b)] There is an open subgroup of finite index $N\leq G$ such that $V=V^N$.
    \end{enumerate}
    Representations with finite action are finite-dimensional, but the converse is false
    in general. However, if $\Gamma$ is a cocompact lattice in $G$,
    then irreducible representations that are weakly contained in  $\LL{\Gamma\leftmod G}$
    have finite action if and only if they are finite dimensional \cite[Lm.~6.20]{Fi15A}.

\subsection{A Criterion for Being Ramanujan}
\label{subsec:being-Ramanujan}

	
	 \begin{thm}[{\cite[Th.~6.22]{Fi15A}}]\label{TH:Ramanujan-criterion}
        Write $F=\llf\circ S$ where $S$ is as in Definition~\ref{DF:elementary-functor}.
        Let $x_1,x_2,\dots,x_t$ be  representatives of the $G$-orbits in $S \calX$,
        and let $K_n=\Stab_G(x_n)$ ($1\leq n\leq t$). Then:
        \begin{enumerate}
            \item[(i)] $\Gamma\leftmod\calX$ is $F$-Ramanujan
            if and only if
            every irreducible unitary representation  $V\wc \LL{\Gamma\leftmod G}$
            satisfying $V^{K_1}+\dots+V^{K_t}\neq 0$
            is tempered (i.e.\ $V\wc \LL{1\leftmod G}$) or has finite action.
            \item[(ii)] $\Gamma\leftmod\calX$ is completely Ramanujan if and only if
            every irreducible unitary representation  $V\wc \LL{\Gamma\leftmod G}$
            is tempered or has finite action.
        \end{enumerate}
        When $\Gamma\leftmod \calX$ is finite, one can replace ``finite action'' with ``finite dimension''
        and ``$V\wc \LL{\Gamma\leftmod G}$'' with ``$V\leq \LL{\Gamma\leftmod G}$''.
    \end{thm}
    
    The criterion is a consequence of the following theorem, which gives more information
    about how the $F$-spectrum of $\Gamma\leftmod\calX$ is related to the irreducible
    subrepresentations of $\LL{\Gamma\leftmod \calX}$.
	
	\begin{thm}[{\cite[Th.~6.21]{Fi15A}}]\label{TH:spectrum-correspondence}
        Keep the notation of Theorem~\ref{TH:Ramanujan-criterion}.
        For any subset $T$ of the unitary dual
        $\udual{G}$,
        define
        \[T^{(K_1,\dots,K_t)}=\{[V]\in T\suchthat V^{K_1}+\dots+V^{K_t}\neq 0\}. \]
        Then there exists an additive functor
		\[
		\calF:\Rep[u]{G}\to \Rep[u]{\Alg{\catC,F}{}}
		\]
        with the following properties:
        \begin{enumerate}
            \item[(i)]
				$\calF$ induces an embedding
				$[V]\mapsto [\calF V]:\udual{G}^{(K_1,\dots,K_t)}\to \tAlg{\catC,F}{}$,
				denoted $\what{\calF}$.
			\item[(ii)]
				For all $[V]\in\udual{G}^{(K_1,\dots,K_t)}$ and
				$V'\in\Rep[u]{G}$, we have $V\wc V'$ $\iff$ $\calF V\wc \calF V'$.
			\item[(iii)]
				$\calF(\LL{\Gamma\leftmod G})=\overline{F(\Gamma\leftmod\calX)}$
				for all $\Gamma\leq_{\calX}G$.
            \item[(iv)] Let $[V]\in\udual{G}^{(K_1,\dots,K_t)}$.
            	Then $\what{F}[V]$ is in $\frakT_{\Alg{\catC,F}{}}$
            	(the trivial $\Alg{\catC,F}{}$-spectrum) 
            	if and only if  $V$ has  finite action.
            	More precisely, we have $\what{\calF}(\Spec_{G}(\llf(N\leftmod G))^{(K_1,\dots,K_t)})=\frakT_{\Alg{\catC,F}{},N}$
            	for any open finite-index subgroup $N\leq G$.
        \end{enumerate}
    \end{thm}
    


	It is also possible to describe $\Alg{\catC,F}{}$ using the \emph{Hecke algebra} of $G$. See
	\cite[Th.~6.10]{Fi15A} for details.
	
\medskip

	We now give some consequences of Theorems~\ref{TH:Ramanujan-criterion}
	and~\ref{TH:spectrum-correspondence}.
	
\medskip

	Marcus, Spielman and Srivastava \cite{MarSpiSri14}
	proved that every bipartite $k$-regular graph $X$ has a Ramanujan double
	cover $X'\to X$. That is, all eigenvalues of the vertex adjacency
	operator of $X'$ not arising from $X$ are in the interval
	$[-2\sqrt{k-1},2\sqrt{k-1}]$. This was extended to covers
	of any prescribed degree by Hall, Puder and Sawin \cite{HallPudSaw16}.
	Using Theorem~\ref{TH:spectrum-correspondence}, one can show the following
	result, which is reminiscent of the Ramanjuan--Petersson conjecture.
	
	\begin{cor}[{\cite[Cor.~6.28, Pr.~6.30]{Fi15A}}]\label{CR:ramanujan-subgroups-of-tree}
		Let $\calX$ be a $k$-regular tree, let $G=\Aut(\calX)$
		and let $H$ be the index-$2$ subgroup of $G$ consisting of
        automorphisms preserving the canonical $2$-coloring of $\calX^{(0)}$.
        Then for any cocompact $\Gamma'\leq H$ and $r\in \N$, there exists a sublattice $\Gamma'\leq \Gamma$ of index $r$ such that
        every irreducible unitary subrepresentation $V$ of the orthogonal complement of
        $\LL{\Gamma\leftmod G}$ in $\LL{\Gamma'\leftmod G}$
        is tempered.
    \end{cor}
    
    The following proposition follows from Theorem~\ref{TH:Ramanujan-criterion} 
    and representation-theoretic properties
    of the relevant group $G$.
    
    \begin{prp}[{\cite[Pr.~6.30, Pr.~6.31]{Fi15A}}]\label{PR:Raman-graph-equiv-conds}
    	(i) Let $\calX$ be a $k$-regular tree, and let $G=\Aut(\calX)$.
    	Then a $k$-regular graph, viewed as a $G$-quotient of $\calX$,
    	is Ramanujan in  dimension $0$ (i.e.\ Ramanujan in the classical sense)
        if and only if it is completely Ramanujan.
        
        (ii) 
        Let $F$ be a local non-archimedean field,
        let $G=\nPGL{F}{d}$ with  $d\in\{2,3\}$, let $\calX$ be the affine Bruhat--Tits building of $G$,
        and let $\Gamma\leq_{\calX}G$.
        Then $\Gamma\leftmod \calX$ is Ramanujan in dimension $0$
        if and only if $\Gamma\leftmod\calX$ is Ramanujan in all dimensions.
    \end{prp}
	

	Let $H$ be an open finite index subgroup of $G$.
	Then $\calX$ is  an almost transitive $H$-complex.
	For $\Gamma\leq_\calX H$, we can consider
	$\Gamma\leftmod \calX$
    both as a $G$-quotient and as an $H$-quotient of $\calX$.
    However, using Theorem~\ref{TH:Ramanujan-criterion}(ii), one can show
    that this does not affect the completely Ramanujan property.
    
    \begin{thm}[{\cite[Th.~6.36]{Fi15A}}]\label{TH:finite-index-subgroup}
        Let $\Gamma\leq_{\calX}H$. Then $\Gamma\leftmod \calX$ is 
        completely Ramanujan as an $H$-quotient of $\calX$
        if and only if $\Gamma\leftmod\calX$ is completely Ramanujan as a $G$-quotient of $\calX$.
    \end{thm}

\section{Ramanujan Complexes}
\label{sec:ram-comp}

	Let $F$ be a local non-archimdean field of positive characteristic,
	let $D$ be a central division $F$-algebra, let $d\geq 2$,
	let $G=\nPGL{D}{d}:=\nGL{D}{d}/\units{F}$,
	and let $\calB_d(D)$ be the affine Bruhat--Tits building of $G$ (we recall its construction below).
	
	Theorem~\ref{TH:Ramanujan-criterion} can be applied
	together with deep results about automorphic representations
	to show that $\calB_d(D)$ has infinitely many non-isomorphic $G$-quotients
	which are completely Ramanujan. This was shown in the case $D=F$
	by Lubotzky, Samuels and Vishne \cite{LubSamVi05}, with the difference
	that they proved that the $G$-quotients were Ramanujan in dimension $0$.
	
\subsection{The Building of $\nPGL{D}{d}$}
\label{subsec:building}

	Recall that $r:=[D:F]^{1/2}$ is an integer called the \emph{degree} of $D$.
	If $\quo{F}$ is an algebraic closure of $F$, then $\nMat{D}{d}\otimes_F\quo{F}\cong \nMat{\quo{F}}{rd}$
	(as $\quo{F}$-algebras). The restriction of the determinant map $\det:\nMat{\quo{F}}{rd}\to\units{\quo{F}}$
	to $\nMat{D}{d}$ is called the \emph{reduced norm} and denoted $\Nrd_{\nMat{D}{d}/F}$.
	It well-defined and takes values in $F$.
	
	Let $\eta:F\onto \Z\cup\{\infty\}$ denote the additive valuation of $F$.
	By \cite[\S12]{MaximalOrders}, $\eta$ extends uniquely
    to an additive valuation  $\eta_D:D\to \R\cup\{\infty\}$ given by:
    \[
    \eta_D(x)=r^{-1}\eta(\Nrd_{D/F}(x))\ .
    \]
    Since the cardinality $q$ of the residue field of $(F,\eta)$ is finite, $\im(\eta_D)=\frac{1}{r}\Z$ and the residue
    division ring of $(D,\eta_D)$ is the Galois field of cardinality $q^r$ \cite[Th.~14.3]{MaximalOrders}. 
    Fix an element $\pi_D\in D$ with $\eta_D(\pi_D)=\frac{1}{r}$ and write
    \[\calO_D=\{x\in D\suchthat \eta_D(x)\geq 0\}\ .\] 
    
    We make $\nGL{D}{d}$ into a topological group by giving it the subspace topology
    induced from the inclusion $\nGL{D}{d}\subseteq \nMat{D}{d}\cong F^{r^2d^2}$,
    and
    give $G=\nPGL{D}{d}:=\nGL{D}{d}/\units{F}$ the quotient topology.
    Then $G$ is an $\ell$-group.
    The quotient map $\nGL{D}{d}\to \nPGL{D}{d}$ is denoted by $g\mapsto \quo{g}$.
    
\medskip

	The affine Bruhat--Tits building of $G$, denoted
	$\calB_d(D)$, is a simplicial complex constructed as follows:
	Let $K$ be the subgroup of $G$ generated by the images of images of $\nGL{\calO_D}{d}$ and
    \[
    \left[\DDotsArr{\pi_D}{\pi_D}\right]
    \]
    in $G$. The vertices of $\calB_d(D)$ are $G/K$.
    To define the edges of $\calB_d(F)$, let
        \[
        g_1=\left[\begin{smallmatrix} \pi_D & & & & \\ & 1 & & & \\ & & 1 & & \\ & & & \ddots & \\ & & & & 1\end{smallmatrix}\right],\quad
        ~ \dots~,\quad
        g_{d-1}=\left[\begin{smallmatrix} \pi_D & & & & \\ & \pi_D & & & \\ & & \ddots & & \\ & & & \pi_D & \\ & & & & 1\end{smallmatrix}\right]
        \in\nGL{F}{d}
        \]
    Two vertices ${g}K,{g'}K\in G/K$ are adjacent if
    \[
    {g^{-1}g'}\in K\cup K\quo{g_1}K\cup K\quo{g_2}K\cup\dots K\quo{g_{d-1}}K\ ,
    \]
    and the $i$-dimensional cells of $\calB_d(D)$ are the $(i+1)$-cliques, namely,
    they are
    sets $\{h_0K,\dots,h_{i+1}K\}\subseteq G/K$
    consisting of pairwise adjacent vertices.
    The resulting complex is indeed a {pure} $(d-1)$-dimensional contractible simplicial complex,
    which carries additional structure making it into an \emph{affine building}; 
    see \cite[\S6.9]{Buildings08AbramBrown} or \cite{AbramNebe02}
    for further details.
    
    There is an obvious left action of $G$ on $\calB_d(D)$, making the latter  into an
    almost transitive $G$-complex.
    
    \begin{example}
            The complex $\calB_2(D)$ is a $q^r+1$ regular tree. 
            Its $G$-quotients  are therefore
            $(q^r+1)$-regular graphs.
    \end{example}

	We now explain why the $0$-dimensional spectrum of $G$-quotients of $\calB_d(F)$
	in the sense of \ref{subsec:spectrum} is equivalent to the spectrum 
	defined by Lubotzky, Samuels and Vishne in \cite{LubSamVi05}
	(basing on Cartwright, Sol\'{e} and \.{Z}uk \cite{CarSolZuk03}).

	The spectrum of Lubotzky, Samuels and Vishne is defined as follows: Consider the map $c:\nGL{D}{d}\to \Z/d\Z$ given by
    \[c(g)={ \eta(\Nrd_{\nMat{D}{d}/F} (g))}+d\Z\ .\]
	It induces a $(d-1)$-coloring  of the \emph{directed} edges of $\calB_d(D)$ given by
    \[
        C_1(\quo{g}K,\quo{g'}K):=c(g^{-1}g')\in \Z/d\Z,\qquad \forall\,g,g'\in\nGL{D}{d}\ .
    \]
    This is a $(d-1)$-coloring because $\im(C_1)=\Z/d\Z-\{0\}$.
    
    The coloring $C_1$ descends to the $G$-quotients of $\calB_d(F)$.
    If $X$ is a such a quotient, we define linear operators
    \[
    a_{1,X},\dots,a_{d-1,X}:\Omega_0^+(X)\to \Omega_0^+(X)
    \]
    by
    	\[
        (a_{i,X}\vphi)x=\sum_{\substack{y\in \vrt{X}\\C_1(x,y)=i}}\vphi y\qquad\forall \,\vphi\in \LL{\vrt{X}},\, x\in \vrt{X}\ .
        \]
        The operators $a_{1,X},\dots,a_{d-1,X}$ are called the \emph{colored adjacency operators} or \emph{Hecke operators}
        of $X$.
        It turns out that they commute with each other and that $a_{i,X}^*=a_{d-i,X}$ for all $i$.
        According to Lubotzky, Samuels and Vishne \cite{LubSamVi05}, 
        when $D=F$,
        the spectrum of a $G$-quotient $X$ of $\calB_d(D)$ is 
        \[
        \Spec(a_{1,X},\dots,a_{d-1,X})\subseteq \C^{d-1}\ .
        \]
        This definition also makes sense when $D\neq F$.

        Now, write $\catC=\catC(G,\calB_d(D))$ and  $a_i=\{a_{i,X}\}_{X\in\catC}$.
        Then $a_1,\dots,a_{d-1}$ are operators associated with $(\catC,\Omega_0^+)$
        and they span a commutative $*$-subalgebra  of $A(\catC,\Omega_0^+)$.
      	When $D=F$, it is known  that $\C[a_1,\dots,a_{d-1}]=A(\catC,\Omega_0^+)$ 
        (\cite[Ch.~V]{Macdonald95}; see also \cite[Ex.~6.14]{Fi15A}).
        Thus, by virtue of Proposition~\ref{PR:unitary-dual-topology},
        the spectrum of Lubotzky, Samuels and Vishne \cite{LubSamVi05}
        is equivalent to the $0$-dimensional spectrum of $G$-quotients of $\calB_d(F)$.

\medskip
        
        Lubotzky, Samuels and Vishne \cite{LubSamVi05} also defined Ramanujan $G$-quotients
        of $\calB_d(F)$ as those quotients whose spectrum is contained
        in union of the spectrum of $\calB_d(F)$ (which they determined
        in \cite[Th.~2.11]{LubSamVi05}) with a certain set of points in $\C^{d-1}$
        called the \emph{trivial spectrum} (see \cite[\S2.3, Def.~1.1]{LubSamVi05}).
        The $G$-quotients which are Ramanujan in this sense are precisely
        the $G$-quotients which are Ramanujan in dimension $0$ according to our definition
        \cite[Ex.~5.15(ii)]{Fi15A}.
        
        \begin{remark}
        	The assertion $\C[a_1,\dots,a_{d-1}]=A(\catC,\Omega_0^+)$ seems
        	to be correct for general $D$. We were unable to find a source, however.
        \end{remark}
        
\subsection{Ramanujan Quotients}
\label{subsec:ramanujan-exist}

		Keep the notation of \ref{subsec:building}.
		We now  state our main result, which gives particular infinite families
		of completely 
		Ramanujan $G$-quotients of $\calB_d(D)$. We remind the reader
		that we assume  $\Char F>0$.
		
\medskip

		We introduce additional notation: There is a global field
		$k$ with a place $\eta$ such that the completion of $k$
		at $\eta$, denoted $k_\eta$, is $F$.
		Let $\calV$ be the set of places of $k$. The additive
		valuation corresponding to $\nu\in\calV$ is also denoted $\nu$.
		We further write
		\begin{align*}
		\calO_\nu &= \{\alpha\in k_\nu\suchthat \nu(\alpha)\geq 0\}, \\
		R&=\{\alpha\in k\suchthat \text{$\nu(\alpha)\geq 0$ for all $\eta\neq \nu\in \calV$}\}\ .
		\end{align*}
		
		Let $E$
		be a central division $k$-algebra of dimension $(rd)^2$.
		Then for every $\nu\in\calV$, there is a central
		division $k_\nu$-algebra $D_\nu$ and $m_\nu\in \N$
		such that $E_\nu:=E\otimes_kk_\nu\cong \nMat{D_\nu}{m_\nu}$.
		Suppose that $E$ is chosen such that
		\begin{itemize}
        \item $E_\eta\cong \nMat{D}{d}$, i.e.\ $D_\eta\cong D$ and $m_\eta=d$.
        \item There is $\theta\in \calV$ such that $E_\theta$ is a division ring, i.e.\ $m_\theta=1$.
    	\end{itemize}
    	Existence of a suitable $E$ for any prescribed $D$ and $d$ follows from the
		Albert--Brauer--Hasse--Noether Theorem (\cite[Rem.~32.12(ii)]{MaximalOrders}, for instance).

		The functor $A\mapsto \Aut_{A\textrm{-}\mathrm{alg}}(E\otimes_kA)$
		from commutative $k$-algebras to groups is representable by an affine group scheme
		over $k$, denoted $\uPGL_{1,E}$. We write $\bfH=\uPGL_{1,E}$ for brevity.
		By the Skolem--Noether Theorem,
		$\bfH(L)=\units{(E\otimes_kL)}/\units{L}$ for every field extension $L/k$.
		Fix a closed
		embedding $j:\bfH\to \uGL_{n}$.
		For any integral domain $S$ whose fraction field $L$ contains $k$,
		and for any $I\idealof S$, we write
		\begin{align*}
		\bfH(S)&=j(\bfH(L))\cap \uGL_n(S)\\
		\bfH(S,I)&=\ker(\bfH(S)\to \uGL_n(S)\to \uGL_n(S/I))\ .
		\end{align*}
		
		We assume that
		\begin{itemize}
			\item $\bfH(\calO_\theta)=\im(\units{\calO_{E_\theta}}\to \units{E_\theta}/\units{k_\theta}=\bfH(k_\theta))$
		\end{itemize}
		where $\calO_{E_\theta}$ is defined as in \ref{subsec:building}.
		The existence of an embedding $j:\bfH\to \uGL_n$ with this property
		is shown in \cite[\S5]{LubSamVi05}.
		
		Finally, recall that the ideals of $R$ correspond to functions
		$\vec{n}: \calV-\{\eta\}\to\N\cup\{0\}$ of finite support. The ideal
		corresponding to $\vec{n}$ is
		\[
		I(\vec{n})=
		\{\alpha\in R\suchthat \text{$\nu(\alpha)\geq \vec{n}(\alpha)$ for all $\nu\in\calV-\{\eta\}$}\}\ .
		\]
		We write
		\[
		\Gamma(\vec{n})=\bfH(R,I(\vec{n}))\ .
		\]
		Since $\bfH$ is $k$-anisotropic, $\Gamma(\vec{n})$ is a cocompact lattice in 
		$\bfH(k_\eta)=\nPGL{D}{d}$ (see \cite[\S{}7D]{Fi15A}).
		For every $\nu\in\calV-\{\eta\}$, there
		is $n_0\in\N$ such that $\Gamma(\vec{n})\leq_{\calB_d(D)}\nPGL{D}{d}$
		whenever $\vec{n}(\nu)\geq n_0$ \cite[Rm.~7.26(ii)]{Fi15A}.

		\begin{thm}[{\cite[Th.~7.22]{Fi15A}}]\label{TH:ram-quo-exists}
			Let $\vec{n}:\calV-\{\eta\}\to\N\cup\{0\}$ be a function
			of finite support such that $\Gamma(\vec{n})\leq_{\calB_d(D)}\nPGL{D}{d}$.
			Assume that either
			\begin{enumerate}
			\item[(1)] $\vec{n}(\theta)=0$, or 
			\item[(2)] $D=F$ and $d$ is prime.
			\end{enumerate}
			Then $\Gamma(\vec{n})\leftmod\calB_d(D)$ is completely Ramanujan.
		\end{thm}
		
		When $D=F$, Theorem~\ref{TH:ram-quo-exists}  
		is just Theorem~1.2 of \cite{LubSamVi05} with the difference
		that we show complete Ramanujan-ness whereas  \cite{LubSamVi05} shows Ramanujan-ness in dimension $0$
		(cf.\ \ref{subsec:building}).

		\begin{example}[{\cite[Ex.~7.23]{Fi15A}}]
        	For every $\Gamma=\Gamma(\vec{n})$ as in Theorem~\ref{TH:ram-quo-exists}, the spectrum
        	of the $i$-dimensional Laplacian $\Delta_i$ of $\Gamma\leftmod\calB_d(D)$ is contained
        	in the union of the spectrum of the $i$-dimensional Laplacian of $\calB_d(D)$  
        	with the trivial spectrum of $\Delta_i$ (which is $\frakT_{\C[\Delta_i]}$,
        	where $\C[\Delta_i]$ is the subalgebra of $\Alg{\catC(\nPGL{D}{d},\calB_d(D)),\Omega_i^-}{}$
        	spanned by $\Delta_i$; use  Proposition~\ref{PR:unitary-dual-topology}
        	and Proposition~\ref{PR:Ramanujan-well-behaved}).
        	This also holds for the adjacency operators considered in Example~\ref{EX:higher-dimensional-adj-ops}.
        	Therefore,
        	when $d=2$, the quotient $\Gamma\leftmod \calB_2(D)$ is a Ramanujan $(q^r+1)$-regular
        	graph.
    	\end{example}
    	
    	Theorem~\ref{TH:ram-quo-exists} is proved
    	by using the criterion of Theorem~\ref{TH:Ramanujan-criterion}(ii)
    	together with deep results about automorphic representations, particularly
    	the proof of the Rama\-nujan--Peters\-son conjecture for $\uGL_n$
    	in positive characteristic due to Lafforgue \cite{Laff02},
    	and the global Jacquet--Langlands correspondence for $\uGL_n$
    	in positive characteristic, established in \cite{BadulRoch14}.
    	
    	In the next subsection, we give brief details about how to translate the statement of 
    	Theorem~\ref{TH:ram-quo-exists} into a statement about automorphic representations, which
    	can then be proved using  results from \cite{Laff02}, \cite{BadulRoch14} and related
    	works. The details of the latter can be found in \cite[\S{}7F--\S{}7I]{Fi15A}.
    	
\subsection{Automorphic Representations}

	Keeping all previous notation, let $S\subseteq \calV$ be a finite
	subset and let $\bbA^S$ denote the $k$-adeles away from $S$, namely
	\[
    \bbA^S = \prod'_{\nu\in \calV-S}k_\nu :=\left\{(a_\nu)_{\nu}\in\prod_{\nu\in\calV-S}k_\nu\suchthat \text{$a_\nu\in \calO_\nu$ for almost
    all $\nu$}\right\}\ .
    \]
    We give
    $\prod_{\nu\in\calV-S}\calO_\nu$ the product topology, and topologize $\bbA^S$ by viewing it
    as a disjoint union of (additive) cosets of $\prod_{\nu\in\calV}\calO_\nu$.
    If $\bfG$ is an algebraic group over $k$, then we topologize $\bfG(\bbA^S)$
	by choosing a closed embedding $\bfG\to \uSL_n$ and giving $\bfG(\bbA^S)$
	the topology induced from $\uSL_n(\bbA)\subseteq \nMat{\bbA}{n}\cong \bbA^{n^2}$.
	This makes $\bfG(\bbA^S)$ into an $\ell$-group\footnote{
		When $\Char k=0$ and $S$ does not contain all archimedean places, 
		$\bfG(\bbA^S)$ is not an $\ell$-group but rather a locally compact group. 
		Nevertheless, the discussion to follow applies when $\Char k=0$ after some modifications. 	
	}; its topology is independent of the embedding.
    
    As usual, $k$ is embedded diagonally in $\bbA^S$, and we write $\bbA:=\bbA^{\emptyset}$.
    We shall occasionally view $\bbA$ as $k_\eta\times \bbA^{\{\eta\}}$.

\medskip

	Let $\bfG$ be a semisimple algebraic group over $k$ with 
	center $\bfZ$. Recall that an automorphic
	representation of $\bfG$ is an irreducible representation 
	$V$ of $\bfG(\bbA)$ that is weakly contained in 
	$\LL{\bfG(k)\leftmod\bfG(\bbA)}$ (cf.\ Example~\ref{EX:right-reg-rep}).
	Every such representation can be written as a restricted tensor product
	$V=\bigotimes'_{\nu\in V_\nu}V_\nu$ with $V_\nu\in\Irr[u]{\bfG(k_\nu)}$; see 
	\cite{Flath79} or \cite[\S{}7C]{Fi15A}. The factor $V_\nu$ is called the 
	\emph{$\nu$-local factor} of $V$.
	

\medskip
 
    Assume $\bfG$ is $k$-anisotropic, and
    choose a compact open subgroup $K^\eta\leq \bfG(\bbA^{\{\eta\}})$. 
    Our assumption implies that 
    $\bfG(k)\leftmod\bfG(\bbA)$ is compact (\cite[Th.~5.5]{PlatRapin94}, \cite[Cr.~2.2.7]{Harder69}),
    so the double coset space
	\[\bfG(k)\leftmod\bfG(\bbA)/(\bfG(k_\eta)\times K^{\eta})\]
	is compact and discrete, hence finite. 
	Let $(1,g_1),\dots,(1,g_t)\in \bfG(k_\eta)\times \bfG(\bbA^{\{\eta\}})$ be representatives
	for the double cosets.
	For each $1\leq i\leq t$, define
	\[
	\Gamma_i=\bfG(k)\cap (\bfG(k_\eta)\times g_i K^{\eta}g_i^{-1})
	\]
	and view $\Gamma_i$ as a subgroup of $\bfG(k_\eta)$. It is a standard fact that there is an
	isomorphism of topological (right) $\bfG(k_\eta)$-spaces
	\begin{equation}\label{EQ:Gamma-decomposition}
	\bigsqcup_{i=1}^t\Gamma_i\leftmod \bfG(k_\eta)\to \bfG(k)\leftmod\bfG(\bbA)/(1\times K^{\eta})
	\end{equation}
	given by sending $\Gamma_i g$ to $\bfG(k)(g,g_i)(1\times K^{\eta})$.
	In particular, $\Gamma_i$ is a cocompact lattice in $\bfG(k_\eta)$ for all $1\leq i\leq t$.

	Assume further  that $\bfG$ is almost simple, let $G=\bfG(k_\eta)/\bfZ(k_\eta)$ and write
    \[
    \quo{\Gamma}_i=\im(\Gamma_i\to G)\ .
    \]
    Then $G$ acts faithfully on the affine Bruhat--Tits building $\calB$ of $\bfG(k_\eta)$, making it into
    an almost transitive $G$-complex (Example~\ref{EX:G-complex-building-III}). 
    Using \eqref{EQ:Gamma-decomposition} and Theorem~\ref{TH:Ramanujan-criterion}, one can show:
    
    \begin{thm}[{\cite[Th.~7.4]{Fi15A}}]\label{TH:automorphic-ramanujan}
        Let $F:\catC(G,\calB)\to \catPHil$ be an elementary functor (e.g.\
        $\Omega_i^+$ or $\Omega_i^\pm$), write $F\cong \llf \circ S$
        as in Definition~\ref{DF:elementary-functor},  let $x_1,\dots,x_s$ be representatives
        for the $G$-orbits in $S\calB$, and let $L_j=\Stab_{\bfG(k_\eta)}(x_j)\times K^\eta$ ($1\leq j\leq s$).
        Assume that for any automorphic representation $V=\bigotimes'_{\nu}V_\nu$
        of $\bfG$
        with $V^{L_1}+\dots V^{L_s}\neq 0$
        (resp.\ $V^{\bfZ(k_\eta)\times K^\eta}\neq 0$), the local factor $V_\eta$ is tempered or finite-dimensional.
        Then $\quo{\Gamma}_i\leftmod \calB$ is $F$-Ramanujan
        (resp.\ completely Ramanujan)
        for every $1\leq i\leq t$ such that $\quo{\Gamma}_i\leq_{\calB} G$.
        The converse holds when $\quo{\Gamma}_i\leq_\calB G$ for all $1\leq i\leq t$.
    \end{thm}
    
    Applying Theorem~\ref{TH:automorphic-ramanujan} with $\bfG=\bfH$ (cf.~\ref{subsec:ramanujan-exist}), 
    a suitable
    $K^\eta$ and $g_1=1$ allows one to translate the statement of Theorem~\ref{TH:ram-quo-exists}
    into a statement about the automorphic representations of $\bfH$, which can be proved
    using powerful results about the latter (see \cite[\S{}7F--\S{}7I]{Fi15A}).
    
    We mention here several places where such ideas were applied in the literature, sometimes implicitly
    or in an equivalent formulation:
    \begin{itemize}
        \item Lubotzkly, Phillips and Sarnak  \cite{LubPhiSar88}, and independently
        Margulis \cite{Marg88}, constructed infinite families
        of Ramanujan $(p+1)$-regular graphs for every prime $p$ using results of Eichler \cite{Eichler54} and Igusa \cite{Igusa59}
        about modular forms. (See also Delinge's proof of the Ramanujan--Petersson conjecture
        for modular forms \cite{Deligne74}.)
        In the previous setting, this corresponds to taking $k=\Q$
        and $\bfG$ to be an inner form of $\bPGL_2$ which  splits over $k_\eta$.
        \item Morgenstern \cite{Morg94} used Drinfeld's proof of the Ramanujan--Petersson conjecture
        for $\bGL_2$ when $\Char k>0$ \cite{Drinfel88}
        to construct  infinite families of Ramanujan ${(q+1)}$-regular graphs for every prime power $q$. Again, the
        corresponding group $\bfG$ is an inner form of $\bPGL_2$.
        \item Lubotzky, Samuels and Vishne \cite{LubSamVi05} applied Lafforgue's proof of the Ra\-ma\-nu\-jan--Petersson conjecture
        for $\bGL_d$ when $\Char k>0$ \cite{Laff02} to construct infinite families of quotients
        of $\calB_d(F)$ which are Ramanujan in dimension $0$. The corresponding group $\bfG$ is an inner form of $\bPGL_n$
        which splits over $k_\eta$.
        \change{4}{
        \item Li \cite{Li04} independently gave similar constructions of Ramanujan complexes, using
        results of Laumon, Rapoport and Stuhler, who proved a special case of the Ramanujan--Petersson
        conjecture for anisotropic inner forms of $\uGL_n$  \cite[Th.~14.12]{LaRaSt93}.
        }
        \item Ballantine and Ciubotaru  \cite{BallCiub11} constructed infinite families of Ramanujan $(q+1,q^3+1)$-biregular graphs
        for every prime power $q$. The corresponding group $\bfG$ is an inner form of $\mathbf{SU}(3)$, and they
        use the classification of the automorphic spectrum of $\bfG$ due to
        Rogawski \cite{Roga90}.
    \end{itemize}
    We hope  our work will facilitate further results of this kind.
	
\bibliographystyle{plain}
\bibliography{MyBib_16_05}

\end{document}